\documentclass[11pt,reqno]{amsart}

\usepackage[utf8]{inputenc}
\usepackage[margin=1.25in]{geometry}
\usepackage{hyperref}
\usepackage{appendix}
\usepackage{amsfonts}
\usepackage{amsthm}
\usepackage{amssymb}
\usepackage{stmaryrd} 
\usepackage{amsmath}
\usepackage{amsthm}
\usepackage[dvipsnames]{xcolor}
\usepackage{mathrsfs}
\usepackage{lipsum} 

\theoremstyle{plain}
\newtheorem{theorem}{Theorem}[section]
\newtheorem{corollary}[theorem]{Corollary}
\newtheorem{lemma}[theorem]{Lemma}

\newtheorem{Assumption}[theorem]{Assumption}

\newtheorem{Definition}[theorem]{Definition}

\newtheorem{fact}[theorem]{Fact}

\newtheorem{notation}[theorem]{Notation}
\newtheorem{Claim}[theorem]{Claim}

\theoremstyle{remark}
\newtheorem{example}[theorem]{Example}

\newtheorem{remark}[theorem]{Remark}

\usepackage{biblatex} 
\numberwithin{equation}{section}
\title[Invertibility for random band matrices]{Invertibility for non-Hermitian and symmetric random band matrices with sublinear bandwidth and discrete entries}

\author{Yi HAN}
\address{Department of Mathematics, Massachusetts Institute of Technology, Cambridge, MA
}
\email{hanyi16@mit.edu}

\begin{document}

\begin{abstract}
A well-known result in random matrix theory, proven by Kahn, Komlós and Szemerédi in 1995, states that a square random matrix with i.i.d. uniform $\{\pm 1\}$ entries is invertible with probability $1-\exp(-\Omega(n))$. As a natural generalization of the model, we consider the invertibility of a class of random band matrices with independent entries where the bandwidth $d_n$ scales like $n^\alpha$, for some $\alpha\in(0,1)$. The band matrix model we consider is sufficiently general and covers existing models such as the block band matrix and periodic band matrix, allowing great flexibility in the variance profile. As the bandwidth is sublinear in the dimension, estimating the invertibility and least singular values of these matrices is a well-known open problem. We make progress towards the invertibility problem by showing that, when $\alpha>\frac{2}{3}$ and when the random variables are i.i.d. uniformly distributed on $\{\pm 1+c\}$ for any fixed integer $c$, then the band matrix is invertible with probability $1-\exp(-\Omega(n^{\alpha/2}))$. Previously, even invertibility with probability $1-o(1)$ was not known for these band matrix models except in the very special case of block band matrices. 

We then extend the invertibility result to symmetric random band matrices with integer entries, and prove the same non-singularity 
probability estimate whenever $\alpha>\frac{2}{3}$.

\end{abstract}

\maketitle

\section{Introduction}

In this paper, we investigate the invertibility for a wide class of combinatorial random matrices with an inhomogeneous variance profile. Before stating the main result on band matrices, we begin with a historical account on the invertibility of $n\times n$ random matrices with i.i.d. entries, which is much better understood in the past thirty years. Let $M_n$ be an $n\times n$ matrix with each entry independently and uniformly distributed on $\{-1,1\}$, and denote by $P_n:=\mathbb{P}(\det M_n=0)$. The first proof of asymptotic non-singularity, that $\lim_{n\to\infty}P_n=0$, was given by Komlós \cite{komlos1967determinant}. Then Kahn, Komlós and Szemerédi \cite{kahn1995probability} provided the first exponential bound on $P_n$: $P_n\leq c^n$ for $n$ sufficiently large, and where they can take $c=0.999$. Tao and Vu \cite{tao2007singularity} improved the value of $c$ to be $c=3/4+o(1)$, and Bourgain, Vu, and Wood \cite{bourgain2010singularity} further improved to $c=\sqrt{2}/2+o(1)$. Finally, Tikhomirov \cite{tikhomirov2020singularity} obtained the optimal value $c=\frac{1}{2}+o(1)$. A closely related model is the random symmetric matrix model $M_n^{sym}$ with uniform $\{-1,1\}$ entries, whose invertibility was first studied in \cite{costello2006random}; and Campos, Jenssen, Michelen, and Julian Sahasrabudhe \cite{campos2025singularity} recently proved that the singularity probability of $M_n^{sym}$ is also exponentially small, namely we can find $c'>0$ such that $\mathbb{P}(\det M_n^{sym}=0)\leq\exp(-c'n)$.

In parallel lines, much recent progress has been made for random matrices with an inhomogeneous variance profile, in the dense case where the matrix typically has $\Theta(n^2)$ nonzero entries. See \cite{rudelson2016singular} for least singular value estimates in the Gaussian case and see \cite{MR3878135}, \cite{MR3857860}, \cite{cook2article} for the general case where the circular law was proven when the matrix has a doubly stochastic variance profile. However, the case of sparse structured matrices was not well understood, and proving the convergence of its empirical spectral measure was listed as an open problem in the survey \cite{tikhomirov2022quantitative}.

For sparse structured random matrices, there are essentially only two recent papers \cite{jain2021circular} and \cite{tikhomirov2023pseudospectrum} studying the invertibility and pseudospectrum. The paper \cite{jain2021circular} by Jain, Jana, Luh and O’Rourke considered a very special type of matrix model called block band matrices (see Definition \ref{blockcanonicalform}), and proved convergence of empirical eigenvalue density for bandwidth $d_n\gg n^{33/34}$ utilizing an estimate they developed for the least singular value. Notably, the work \cite{jain2021circular} covers discrete entry laws such as the uniform distribution on $\{\pm 1\}$. However, the special structure of block band matrices is utilized in an essential way, resisting a generalization to more general models such as the periodic band matrices in Definition \ref{periodicbandmatrixform}. Another recent paper by Tikhomirov \cite{tikhomirov2023pseudospectrum}
obtained minimal singular value estimates for (the shifted version of) a wide class of random matrices with general variance profile, but assuming the random variables have a bounded density. Circular law for periodic band matrices with Gaussian entry distribution was proven as a corollary, for bandwidth up to $d_n\gg n^{33/34}$. The range of $d_n$ for circular law to hold was recently significantly extended, see Remark \ref{afterthispaperwas}. However, \cite{tikhomirov2023pseudospectrum} completely left open the case where the nonzero entries are discrete valued, say over $\{\pm 1\}$. And to our best knowledge, even asymptotic nonsingularity for general random band matrices with discrete entry distribution is completely unexplored. (The spectral outliers for random band matrices are much better studied; see \cite{bandeira2024matrix} and \cite{han2024outliers}).

At this point, we should take caution to distinguish our inhomogeneous matrix model with sparse i.i.d. random matrices, the latter model being fairly well-studied. For a sparse i.i.d. matrix of size $n$ and sparsity $p_n\in(0,1)$, the circular law was proven in \cite{rudelson2019sparse} and \cite{sah2025sparse} (see also earlier references therein) whenever $np_n\to\infty$ and invertibility was studied earlier in  \cite{basak2017invertibility}, \cite{basak2021sharp} and for a related model in \cite{ferber2022singularity}. For the case of constant average degree, namely $p_n=\frac{d}{n}$, the limiting spectral measure was recently studied in \cite{sah2023limiting}. However, the techniques in these papers rely crucially on the homogeneous variance structure of the matrices and thus do not generalize transparently to our inhomogeneous setting.

The objective of this paper is to take a step into the invertibility problem for general non-Hermitian random band matrices with discrete entry distribution. We will study a class of square random matrices $A_n=(a_{ij})_{1\leq i,j\leq n}$ with independent entries, and that nonzero random variables are attached along a band close to the main diagonal, with the bandwidth $d_n$ sublinearly growing in $n$. We can also attach arbitrary integer-valued random variables at locations away from this band. Our model covers block band and periodic band matrices, although being slightly less general than the setting of \cite{tikhomirov2023pseudospectrum} where they consider any doubly stochastic variance profile. Random band matrices, at least for the Hermitian version, are good models interpolating the de-localized phase of GOE/GUE ensembles with mean-field interaction, and the localized phase where Anderson localization takes place. The threshold for localization-delocalization transition is conjectured to be $d_n\sim \sqrt{n}$. A very recent breakthrough \cite{yau2025delocalization} confirmed that when $d_n\gg n^{1/2+c}$, the Hermitian random band matrix exhibits complete eigenvalue delocalization, universality of local eigenvalue statistics and quantum unique ergodicity, among other things. The same results are proven for more general band matrix models by Erdős and Riabov \cite{erdHos2025zigzag}. The localized phase where $d_n\ll\sqrt{n}$ is also recently studied in \cite{cipolloni2024dynamical}, \cite{chen2022random}, \cite{goldstein2022fluctuations}, \cite{drogin2025localization}. A more comprehensive survey on prior results in  (Hermitian) random band matrices can be found in \cite{MR3966510} and \cite{MR3966868}. 

Although none of these results on Hermitian band matrices are proven for the non-Hermitian model, we can make the following postulate: invertibility with high probability for a random matrix with only discrete entry law should be a mean field type behavior and relies crucially on the mean field structure. That is, the fact that the square matrix $M_n$ with discrete $\{\pm 1\}$ entry is invertible with very high probability ($P_n\leq c^n$ as cited) stems from the fact that we are in the delocalized phase with strong mean-field interaction. This can be contrasted with random Schrödinger operators in the localized phase (see for example \cite{kirsch2007invitation}), where the Wegner estimate, a quantitative version of invertibility, can generally be established only when the random potentials have a bounded density. Now we turn back to the setting of Theorem \ref{maintheorem}, the discussion of non-Hermitian band matrices. When the band matrix entries have a bounded density, it is possibly not hard to prove that the matrix is almost surely invertible for bandwidth $d_n$ much smaller than $n^{1/2}$: this is because we use the bounded density assumption of each entry, so the system may not be required to lie in the delocalized phase. However, if the band matrix has only discrete entries, and $d_n\ll n^{1/2}$ so the matrix is in the localized regime, we are not sure whether the lack of mean field behavior can still guarantee invertibility with high probability for this band matrix with discrete entry law.

To conclude, from a viewpoint of Anderson localization-delocalization transition, there are two regimes, either (i) $d_n\gg n^{1/2}$, and we are in a delocalized regime, then we can reasonably believe that this band matrix with discrete entry law is invertible with high probability, and the proof uses tools for mean field random matrices; or (ii) in the $d_n\ll n^{1/2}$ scenario, we do not know if invertibility for a discrete entry law persists in general, and even if it does, we feel that the techniques from mean field random matrices will no longer be applicable as the matrix deviates too far from its mean field regime (and other special structures of the matrix will be more useful, see Remark \ref{lines9091}). In this paper we will focus only on scenario (i) where invertibility is heuristically believed to hold. In particular, the main result of this paper is  stated for $d_n\gg n^{1/2}$ (we actually need $d_n\gg n^{2/3}$), and we leave open the possibility of invertibility for $d_n\ll n^{1/2}$.

The main result of this paper is as follows.

\begin{theorem}\label{maintheorem}
    Fix $\alpha\in(0,1)$ and let $d_n$ be a sequence of integers such that $d_n\geq n^\alpha$ when $n$ is sufficiently large. Also assume that $d_n\leq n^{1-\beta}$ for some small $\beta>0$.

    Let $A_n=(a_{ij})_{1\leq i,j\leq n}$ be a family of square random matrices with independent entries. We assume that 
    \begin{enumerate}
        \item When $|i-j|\leq d_n$, then we can write $a_{ij}=u_{ij}+c_{ij}$ where $u_{ij}$ are i.i.d. and uniformly distributed on $\{-1,1\}$, and $c_{ij}$ are arbitrary integer-valued random variables that are independent of $u_{ij}$.
        \item All the other random variables $a_{ij}$ for $|i-j|>d_n$ are either identically zero or are arbitrary integer-valued random variables.
    \end{enumerate} Then whenever $\alpha>\frac{2}{3}$, we have 
    $$
\mathbb{P}(A_n\text{ is singular})\leq \exp(-Cn^{\alpha/2}),
    $$ for some universal constant $C>0$.
\end{theorem}

Observe that by choosing the constants $c_{ij}$ arbitrarily, we can assume that $\{a_{ij}\}_{|i-j|\leq d_n}$ are uniform over $\{0,2\}$, or uniform over $\{-1+c,1+c\}$ for any fixed integer $c$.

The assumption on $A_n$ in Theorem \ref{maintheorem} is sufficiently general and covers the following two special cases as mentioned in the introduction.
\begin{example}(Block band matrix)\label{blockbandmatrixform} Let 
    $A_n$ be an $n\times n$ matrix. We say $A_n$ is sampled from a block band matrix ensemble if it has the form
\begin{equation}\label{blockcanonicalform}
A_n=\begin{pmatrix} D_1&U_2&&&T_m\\T_1&D_2&U_3&&\\&T_2&D_3&\ddots&\\&&\ddots&\ddots&U_m \\ U_1&&&T_{m-1}&D_m \end{pmatrix}
\end{equation}
where the unfilled sites are set zero. The blocks $D_1,U_1,T_1,\cdots,D_m,U_m,T_m$ are independent copies of square random matrices with size $d_n$ and having i.i.d. uniform $\{-1,1\}$ entries.
\end{example}

\begin{example}(Periodic band matrix)\label{periodicbandmatrixform} Let 
    $A_n=\{a_{ij}\}$ be an $n\times n$ square matrix. We say $A_n$ is sampled from a periodic band matrix model with bandwidth $d_n$ if $a_{ij}=0$ for any $d_n<|i-j|<n-d_n$, and $a_{ij}$ has uniform $\{-1,1\}$ distribution if $|i-j|\leq d_n$ or $|i-j|\geq n-d_n$.

We can also consider the (modified) band matrix model where $a_{ij}=0$ when $|i-j|\geq d_n$ and $a_{ij}$ have uniform $\{\pm 1\}$ distributions when $|i-j|< d_n$. 
    
\end{example}

\begin{remark}
    The state of the art for the invertibility probability of non-Hermitian i.i.d. matrices is that the singularity probability of the whole matrix is roughly the same as (or at least has the same exponential order as) the probability that two rows of the matrix are identical. In our model, taking $c_{ij}=1$ for all $|i-j|\leq d_n$ and $a_{ij}=0$ for $|i-j|>d_n$, we find that the probability for two rows of $A_n$ to be the same is at least $2^{-2d_n-1}$ (as now the nonzero entries around the band are uniform on $\{0,2\}$). This probability is exponentially smaller than the upper bound $\exp(-Cn^{\alpha/2})$ in  Theorem \ref{maintheorem}, so we have not proven the same (sharp) heuristic for the singularity probability for our band matrix model here. 
\end{remark}

\begin{remark}
    The assumption in Theorem \ref{maintheorem} that $a_{ij}$ are i.i.d. copies of uniform $\{\pm 1\}$ variables for $|i-j|\leq d_n$ is not significant: via a standard modification we can prove the theorem for the general case where the variables $a_{ij}$ for $|i-j|\leq d_n$ are independent (but not necessarily identically distributed) random variables satisfying that, for some $\beta_0>0$ and for any sufficiently large prime number $p\geq 2$,
    $$\sup_{i,j:|i-j|\leq d_n}
\sup_{x\in\mathbb{Z}}\mathbb{P}(a_{ij}\equiv x\text{ mod }p)\leq 1-\beta_0
.    $$
We do not provide details of this generalization in order to focus on the main ideas of proof.
\end{remark}

Theorem \ref{maintheorem} is an inhomogeneous version (albeit with a weaker probabilistic estimate) of the  Kahn, Komlós and Szemerédi \cite{kahn1995probability} result on invertibility of a square matrix $M_n$ with i.i.d. uniform $\{\pm 1\}$ entries, by showing that invertibility of a square random matrix still holds when we only maintain an $n^{\alpha-1}$ fraction of the total number of random variables in $M_n$, and that this fraction of nonzero random variables are not required to be uniformly distributed among the entries. 
Before this paper, even the qualitative statement of Theorem \ref{maintheorem}, that $A_n$ is nonsingular with probability $1-o(1)$, was not proven for any band matrix model with discrete entries and sublinear bandwidth, with an exception of \cite{jain2021circular}. In \cite{jain2021circular}, Theorem 2.1 they proved that the block band matrix model \eqref{blockcanonicalform} with bandwidth $d_n$ is nonsingular with probability $1-C(d_n)^{-1/2}$. Theorem \ref{maintheorem} improves this singularity probability to exponential type in $n$ whenever $\alpha>\frac{2}{3}$, and yields the first proof of invertibility for periodic band matrices (Example \ref{periodicbandmatrixform}) with $\{\pm 1\}$ entries for any $\alpha\in(\frac{2}{3},1)$.

We believe that proving the invertibility of general band matrices, beyond the special block band matrix setting, requires significant new ideas.
The block band matrix model \eqref{blockcanonicalform} in \cite{jain2021circular} is special in that by the result of \cite{rudelson2008littlewood}, we can assume with high probability that the off-diagonal blocks $T_1,\cdots,T_{m-1}$ and $U_2,\cdots,U_m$ in \eqref{blockcanonicalform} have a bounded inverse, so that we can solve each component of a vector $v$ in the kernel of $A_n$ using iteratively the inverse of these off-diagonal blocks. This special feature is not shared by any other model, such as the periodic band matrix model, whose invertibility is not proven before this work.

The proof of Theorem \ref{maintheorem}
relies on reducing the invertibility problem to a finite field, and using the inverse Littlewood-Offord counting theorems in the finite field setting. Similar ideas have appeared before in \cite{campos2022singularity}, \cite{ferber2021counting}, \cite{campos2021singularity}, but the new feature here is we have to apply the inverse counting theorem for roughly $n/d_n$ times altogether to each component of a vector $v$, and the anti-concentration property for each component of the vector $v$ can be considerably different from one another.

A closely related problem to Theorem \ref{maintheorem} is to estimate the least singular value of $A_n$, and of $zI_n-A_n$ for $z\in\mathbb{C}$ and $I_n$ being the identity matrix. For this purpose, an approach invented by Rudelson and Vershynin \cite{rudelson2008littlewood} is to use the inverse Littlewood-Offord theorem on real numbers instead of the finite field, and use the notion of essential least common denominator (LCD) of a vector invented in \cite{rudelson2008littlewood}. A major problem arises when we apply this approach to band matrices: consider for example the block band matrix model \eqref{blockcanonicalform} (for a general variance profile we actually condition on some other random variables and work only with random variables in this block band structure), consider the components of a unit vector $v\in\mathbb{S}^{n-1}$ with respect to the columns indexed by $D_1,\cdots,D_m$ respectively and denote each individual component by $v_1,\cdots,v_m$. Then the essential LCD of each $v_1,\cdots,v_m$ may vary significantly from one another in the range $[1,\exp(d_n)]$, and the method of \cite{rudelson2008littlewood} requires that we consider an $\epsilon_i$-net for each vector $v_i$ with $\epsilon_i$ depending on the essential LCD of $v_i$. This leads to fundamental problems since, on a row of $A_n$, we are using three different scales of nets to approximate, and thus only the coarsest approximation scale will be used. This leads to significant overcounting when $\epsilon_{i-1},\epsilon_i,\epsilon_{i+1}$ are far apart from each other. Rather, working on a finite field $\mathbb{F}_p$ makes the counting problem much simpler without any need of an approximation, and the finite field version of the inverse Littlewood-Offord counting theorem is far more effective in the case where $p$ is very large but the anti-concentration property of $v_i$ is very weak. This is the key observation that enables our argument to work through. The price we pay for this finite field approach is that we do not get an estimate on the least singular value of $A_n$, which is left for future study.

\begin{remark}\label{afterthispaperwas}
    After this paper was drafted, two recent preprints \cite{han2025circular} and \cite{han2025circular2} proved circular law for non-Hermitian random band matrices with much smaller bandwidth, namely for $\alpha>1/2$ for a class of general models and for any growing bandwidth for a special block band model. These two papers assume a continuous density distribution for individual entries and use this assumption to derive least singular value lower bounds. The techniques used there are thus orthogonal to the techniques in this paper. 
\end{remark}

\begin{remark}\label{lines9091}
We conjecture that an exponential type estimate on the singularity probability of $A_n$ can be proven whenever $\alpha>\frac{1}{2}$ (widely perceived as the whole delocalization regime) or even whenever $\alpha>0$, provided that the inverse Littlewood-Offord counting problem can be significantly refined when $\alpha>\frac{1}{2}$ and provided that special structures of the matrix $A_n$ can be used when considering the entire regime $\alpha>0$. We also conjecture that we can prove $\mathbb{P}(\det A_n=0)=\exp(-\Omega(n^\alpha))$ for $\alpha$ in certain regimes of $(0,1)$.  
\end{remark}

\subsection{Invertibility for random symmetric band matrices}

The proof technique of this paper is robust enough to extend to symmetric random band matrices with integer entries. The following theorem states our main result for the symmetric matrix model.

\begin{theorem}\label{symmetricmaintheorem}
    Fix $\alpha\in(0,1)$ and let $d_n\leq n$ be a sequence of integers such that $d_n\geq n^\alpha$ when $n$ is sufficiently large.

    Let $S_n=(s_{ij})_{1\leq i,j\leq n}$ be a family of square \textbf{symmetric} random matrices such that $s_{ij}=s_{ji}$ and $(s_{ij})_{1\leq i\leq j\leq n}$ are mutually independent random variables. We assume that 
    \begin{enumerate}
        \item When $|i-j|\leq d_n$ or $|n-i+j|\leq d_n$, then for each $1\leq i\leq j\leq n$ we can write $s_{ij}=u_{ij}+c_{ij}$ where $\{u_{ij}\}_{i\leq j}$ are i.i.d. and uniformly distributed on $\{-1,1\}$, and $c_{ij}$ are arbitrary integer-valued random variables that are independent of $u_{ij}$.
        \item All the other random variables $s_{ij}$ where $\min(|i-j|,|n-i+j|)>d_n$   are either identically zero or are arbitrary integer-valued random variables.\end{enumerate}
    Then whenever $\alpha>\frac{2}{3}$, we have 
    $$
\mathbb{P}(S_n\text{ is singular})\leq \exp(-Cn^{\alpha/2}),
    $$ for some universal constant $C>0$.
\end{theorem}

\begin{remark}\label{cyclicwarparound}
    In contrast to Theorem \ref{maintheorem}, in Theorem \ref{symmetricmaintheorem} we take the (similarly standard) assumption that the variance profile wraps around periodically, so that we measure the cycle distance of $i$ and $j$ on $\mathbb{Z}_n$.  
\end{remark}

\begin{remark}
    The very recent bulk universality result \cite{erdHos2025zigzag} for Hermitian band matrices can possibly be applied to show that whenever $\alpha>\frac{1}{2}$ then $S_n$ is non-singular with probability $1-o(1)$. However, this method does not yield the high probability estimate $1-\exp(-Cn^{\alpha/2})$ we derive here for $\alpha>\frac{2}{3}$, and our condition on variance profile in Theorem \ref{symmetricmaintheorem} is not fully comparable with \cite{erdHos2025zigzag}.
\end{remark}

\section{The proof for the independent case}
We first introduce the main technical tool for the proof of Theorem \ref{maintheorem} in Section \ref{section2.1}, and complete the proof of Theorem \ref{maintheorem} in Section \ref{section2.2}.

\subsection{Finite field inverse Littlewood-Offord counting problem}
\label{section2.1}
We shall use one main result from \cite{campos2022singularity} on an inverse Littlewood-Offord counting problem on finite fields. This approach of studying the inverse Littlewood-Offord counting problem on finite fields was first initiated by  \cite{ferber2021counting} and an earlier version of the result from \cite{campos2022singularity} can be found in \cite{campos2021singularity}.

Let $p$ be an odd prime number and $\mathbb{F}_p:=\mathbb{Z}/p\mathbb{Z}$ the finite field with torsion $p$. Fix an integer $d\in\mathbb{N}$. For a vector $v\in\mathbb{F}_p^d$ and $\mu\in[0,1]$ define the random variable $$X_\mu(v)=\epsilon_1v_1+\cdots+\epsilon_dv_d ,$$ where $\epsilon_i\in\{-1,0,1\}$ are i.i.d. random variables with $\mathbb{P}(\epsilon_i=1)=\mathbb{P}(\epsilon_i=-1)=\mu/2$.

We then denote by \begin{equation}\rho_\mu(v)=\max_{x\in\mathbb{F}_p}\mathbb{P}(X_\mu(v)=x),\end{equation} and denote by \begin{equation}\rho(v):=\rho_1(v).\end{equation} Let $|v|$ denote the number of non-zero coordinates of $v$, and for a subset $T\subset[n]$ we denote by $v_T:=(v_i)_{i\in T}$ the restriction of $v$ to $T$.

We shall use the following monotonicity property of $\rho_\mu(v)$  in $\mu$:

\begin{fact}\label{fact2.111}(see \cite{tao2006additive}, Corollary 7.12).
    If $0\leq \mu\leq \mu'\leq 1$ and at least one of $\mu'\leq\frac{1}{2}$ and $\mu\leq\mu'/4$ hold, then\begin{equation}\label{monotonicity}
 \rho_{\mu'}(v)\leq\rho_\mu(v).       
\end{equation}
    
\end{fact}
In this paper, a typical case we shall often use is where $\mu'=1$ and $\mu\in(0,\frac{1}{4}]$. If none of the two conditions on $\mu$ and $\mu'$ hold, the monotonicity can be reversed, see Fact \ref{greenlandfact}. Throughout the entire paper we will only consider values $\mu\in(0,\frac{1}{4}]$, so that we always have for this value of $\mu$ and any $v$ that
$$
\rho(v)\leq\rho_\mu(v).
$$

For $d'\in\mathbb{N}_+$ and a vector $w=(w_1,\cdots,w_{d'})$ we define the neighborhood of $w$ (relative to $\mu$) by 
$$
N_\mu(w):=\{x\in\mathbb{F}_p:\mathbb{P}(X_\mu(w)=x)\geq \frac{1}{2}\mathbb{P}(X_\mu(w)=0)\}.
$$From this definition we can upper bound $N_\mu(w)$ as follows:
\begin{equation}\label{doublecounting}
|N_\mu(w)|\leq\frac{2}{\rho_\mu(w)}.
\end{equation}

This definition of neighborhood is motivated by the fact that, when $\mu\in(0,\frac{1}{2}]$, the random walk $X_\mu(w)$ has the highest probability of taking the value 0 among all the values in $\mathbb{F}_p$ (see \cite{tao2006additive}, Corollary 7.12). Then $N_\mu(w)$ collects the values $x\in\mathbb{F}_p$ to which the random walk $X_\mu(w)$ is at least half as likely to reach compared to the value 0.

The main inverse theorem that we shall use is Theorem \ref{inversetheorem}, taken from \cite{campos2022singularity}, Theorem 2. The original version of Theorem 2 in \cite{campos2022singularity} actually allows for an additional parameter $k\in\mathbb{N}$. Although this parameter $k$ is fundamental for the main result in \cite{campos2022singularity}, it only leads to negligible improvements in our setting, so we only state the theorem for $k=1$.

\begin{theorem}\label{inversetheorem}
    Let $\mu\in(0,\frac{1}{4}]$, $d\in\mathbb{N}$, $p$ a prime number and $v\in\mathbb{F}_p^d$. Write $D=\frac{2}{\mu}\log \rho_\mu(v)^{-1}$. Suppose that $|v|\geq D$ and $\rho_\mu(v)\geq \frac{2}{p}$. Then we can find $T\subset[d]$ with $|T|\leq D$ such that, if we denote $w=v_T$, then $v_i\in N_\mu(w)$ for all but at most $D$ values of $i\in[n]$. Moreover, we have the following cardinality upper bound: $$|N_\mu(w)|\leq \frac{256}{\rho_\mu(v)}.$$
\end{theorem}

We also need the following Lemma \ref{lemmafourier} on the Fourier side, which can be found in the journal version of \cite{campos2022singularity}, Lemma 5. This type of lemma prevails in the literature, starting from \cite{kahn1995probability}, and see also \cite{tao2006additive}, Section 7. 
\begin{lemma}\label{lemmafourier}
    For given $c>0$ and $\mu'\in(0,1]$, we can find $\mu\in(0,\frac{1}{4}]$ and $K>0$ such that for all $d$-dimensional vector $v$ with $\rho_{\mu'}(v)=\Omega(d/p)$ and $|v|\geq K$, we have 
    $$
\rho_{\mu'}(v)\leq c\rho_\mu(v).    $$
\end{lemma}

Combined with the monotonicity statement in Fact \ref{fact2.111}, we can always take a smaller value of $\mu$ at no cost. We provide a proof of Lemma \ref{lemmafourier} for sake of completeness:
\begin{proof}
    By Lemma 7.14 of \cite{tao2006additive}, for any $0<\mu\leq\mu'\leq 1$ with $\mu\leq\frac{1}{4}$, we have
$$
\rho_{\mu'}(v)=O(\sqrt{\frac{\mu}{\mu'}}\rho_\mu(v))+O(\rho_\mu(v)^{\Theta(\mu'/\mu)}).
$$ Then we can take $\mu$ sufficiently small relative to $\mu'$ so that $\sqrt{\frac{\mu}{\mu'}}\leq c/2$.
   Moreover, by \cite{nguyen2015anti}, Theorem 1.5, whenever $K$ is large enough and $|v|\geq K$, for any fixed $\mu\in(0,1),$ we have $\rho_\mu(v)=O_\mu(\max(\frac{1}{p},\frac{1}{K^{1/4}}))$. Then it suffices to take $K$ sufficiently large and we will get $\rho_{\mu'}(v)\leq c\rho_\mu(v)$.  
\end{proof}

\subsection{The main argument for Theorem \ref{maintheorem}}\label{section2.2}

We begin with some simple observations from linear algebra. For any fixed prime number $p$, $A_n$ is singular on $\mathbb{R}$ implies that $A_n$ is singular on the finite field $\mathbb{F}_p$ since $A_n$ has integer coefficients. Then $A_n$ is singular on $\mathbb{F}_p$ implies that there exists some index $\mathbf{I}\in[n]$ such that the $\mathbf{I}$-th row of $A_n$ lies in the $\mathbb{F}_p$-linear span of all the other rows of $A_n$. 
 We will make a careful choice of the prime number $p$, which roughly has order $\exp(n^{\alpha/2})$, and study the singularity probability of $A_n$ on $\mathbb{F}_p$. The latter probability yields an upper bound for the singularity probability of $A_n$ on $\mathbb{R}$.

From this observation, we let $A_n^\mathbf{I}$ denote the matrix obtained from $A_n$ by setting its $\mathbf{I}$-th row to be zero, and denote by $v^\mathbf{I}$ the $\mathbf{I}$-th row of $A_n$. The rationale here is that, if $v^\mathbf{I}$ lies in the span of the rows of $A_n^\mathbf{I}$, then $v^\mathbf{I}$ must be orthogonal to any nonzero vector $v\in\ker (A_n^\mathbf{I})$ (where orthogonality means that $\langle v^\mathbf{I},v\rangle$=0 on $\mathbb{F}_p$). By independence of $v^\mathbf{I}$ and $A_n^\mathbf{I}$, 
we only need to determine certain \textit{arithmetic structure} for any nonzero vectors $v\in\ker A_n^\mathbf{I}\subset\mathbb{F}_p^n$.

In contrast to the full matrix setting (a square matrix with all i.i.d. entries), the main technical difficulty for random band matrices is that the vector $v^\mathbf{I}$ may only have $2d_n\ll n$ coordinates with nonzero random variables, and thus we have to show that any nonzero vector $v\in\ker A_n^\mathbf{I}$ is nonzero when restricted to \textit{these} $2d_n$ coordinates and the restriction of $v$ to these entries has a sufficiently good anti-concentration property. When $d_n=\Theta(n)$, structure theorems can be found in \cite{rudelson2016singular} and \cite{cook2018lower}. But when $d_n$ is sublinear in $n$, no structural results for $v\in\ker A_n^\mathbf{I}$ were previously available for a general variance profile.

Before stating our main structure theorem on vectors in $\ker(A_n^\mathbf{I})$, we introduce some notations on decomposition of the structure of $A_n$. Although we have allowed for a rather arbitrary variance profile for $A_n$ in Theorem \ref{maintheorem}, the first step of the proof is to replace $A_n$ by another matrix with a block band structure as in \eqref{blockcanonicalform}.

We define $e_n=\lfloor d_n/s\rfloor$
(for some $s\in(3,6)$ to be later fixed)\footnote{The value of $s$ is irrelevant to our later proof, but choosing $s$ larger ($s>3$) will be helpful to certain arguments. We also allow a range for $s\in(3,6)$ so that we can always find a $s$ with $n-e_n\lfloor n/e_n\rfloor\geq e_n/2$.} and consider a partition of $[1,n]$ by intervals of length $e_n$:
$$
[1,n]=\cup_{k=1}^{\lfloor n/e_n\rfloor +1}I_k,\quad I_k=\begin{cases} [(k-1)e_n+1,ke_n],\quad \forall 1\leq k\leq \lfloor n/e_n\rfloor \\ [\lfloor n/e_n\rfloor e_n+1,n],\quad k=\lfloor n/e_n\rfloor+1,
\end{cases}$$ where we choose some $s\in(3,6)$ such that $n-e_n\lfloor n/e_n\rfloor\geq e_n/2$. This ensures that the shortest length of the intervals among $I_k$ is at least half as long as the other intervals among $I_k$. We define the following matrix from $A_n=(a_{ij})$ that will be used for comparison: 
\begin{equation}\label{removedsubmatrix}\begin{aligned}
\mathcal{D}_\mathbf{I}&:=\{a_{ij}\mathbf{1}_{i,j\in I_k\text{ or }i\in I_{k-1},j\in I_k\text{ or } i\in I_k,j\in I_{k-1}\text{ for some }k=1,\cdots,\lfloor n/e_n\rfloor+1}\}\cdot \mathbf{1}_{i\neq\mathbf{I}}.\end{aligned}
\end{equation} (where we interpret $I_0=\emptyset$).

In matrix form, we can write $\mathcal{D}_\mathbf{I}$ in the following block form:
\begin{equation}\label{blockformofd}
    \mathcal{D}_\mathbf{I}=\begin{bmatrix}
        D_1&U_2&0&\cdots&0\\T_1&D_2&U_3&\cdots&0\\0&T_2&D_3&\cdots&0\\\ \vdots&\vdots&\vdots&\ddots&U_{\lfloor n/e_n\rfloor+1}\\0&0&0&T_{\lfloor n/e_n\rfloor}&D_{\lfloor n/e_n\rfloor+1}
    \end{bmatrix},
\end{equation} where by the assumption of Theorem \ref{maintheorem}, the entries in the blocks $D_i,U_i,T_i,1\leq i\leq \lfloor n/e_n\rfloor+1$ are independent and uniformly distributed over $\{\pm 1\}+c_{ij}$ for some fixed location-dependent constants $c_{ij}\in\mathbb{Z}$ that do not matter to us, except the $\mathbf{I}$-th row of $\mathcal{D}_\mathbf{I}$ which is identically 0. In other words, $\mathcal{D}_\mathbf{I}$ is a matrix extracted from $A_n$ by zeroing out its $\mathbf{I}$-th row and zeroing out all the entries not in the main diagonal and first off-diagonal blocks labeled by $I_k,k=1,\cdots,\lfloor n/e_n\rfloor+1$.

The main structural theorem on $\ker A_n^\mathbf{I}$ is stated as follows:
\begin{theorem}\label{mainstructure} In the setting of Theorem \ref{maintheorem},
    there exist constants $\rho>0$, $\tau>0$ such that the following holds with $p=\exp(\rho n^{\alpha/2})$: \footnote{By the prime number theorem, there must exist a prime number in $[\exp(\frac{\rho}{2} n^{\alpha/2}),\exp(\rho n^{\alpha/2})]$, and we simply use $p:=\exp(\rho n^{\alpha/2})$ to denote any one of these prime numbers.}

On the finite field $\mathbb{F}_p$, with probability at least $1-\exp(-\Omega(n))$, the following is true for all nonzero vectors $v\in \ker A_n^\mathbf{I}\subset\mathbb{F}_p^n$: if we denote by $v_{I_{n_\mathbf{I}}}$ the restriction of $v$ to the interval $I_{n_\mathbf{I}}$, then we have, with $\rho(v)$  defined in Section \ref{section2.1},
    $$
\rho(v_{I_{n_{\mathbf{I}}}})\leq\exp(-\tau n^{\alpha/2}).
    $$ Moreover, the quantitative estimate is uniform over the choice of $\mathbf{I}\subset[n]$.
\end{theorem}

The proof of Theorem \ref{mainstructure} in fact shows that $\rho(v_{I_k})\leq\exp(-\tau n^{\alpha/2})$ for each $k\leq \lfloor  n/e_n\rfloor+1$, that is, we obtain strong arithmetic properties for each component of $v$ on an interval with sublinear length.

The proof of Theorem \ref{maintheorem} is immediate from Theorem \ref{mainstructure}:

\begin{proof}[\proofname\ of Theorem \ref{maintheorem} assuming Theorem \ref{mainstructure}] By the discussion at the beginning of Section \ref{section2.2}, we have
$$\begin{aligned}&
\mathbb{P}(\det A_n=0)\leq\mathbb{P}(\det A_n=0 \text{ on }\mathbb{F}_p)\leq\sum_{\mathbf{I}=1}^n\mathbb{P}(v^\mathbf{I}\in \operatorname{Span}((v^j)_{j\neq \mathbf{I}})\text{ on }\mathbb{F}_p)\\&\leq\sum_{\mathbf{I}=1}^n\mathbb{P}(\langle v^\mathbf{I},v\rangle=0\text{ for any }v\neq 0:v\in \ker A_n^\mathbf{I}\subset\mathbb{F}_p^n)\\&\leq\sum_{\mathbf{I}=1}^n\min_{0\neq v\in \ker A_n^\mathbf{I}}\rho(v_{I_{n_\mathbf{I}}})\leq n\exp(-\Omega(n^{\alpha/2}))+n\exp(-\Omega(n)),\end{aligned}
$$
    where the last inequality follows from conditioning on the random variables of $v^\mathbf{I}$ indexed by $[n]\setminus I_{n_\mathbf{I}}$ and partitioning the event into two cases: either the conclusion of Theorem \ref{mainstructure} holds for all $v_{I_{n_\mathbf{I}}}$ and all $\mathbf{I}$, or the conclusion of the theorem does not hold for one choice of $\mathbf{I}$, and the latter probability is at most $n\exp(-\Omega(n))$. This concludes the proof.  
\end{proof}

Before outlining the proof of Theorem \ref{mainstructure}, we describe why we assume $\alpha>\frac{2}{3}$:

\begin{remark}
(The assumptions in Theorem \ref{mainstructure}, the main idea and its limitations)

\textbf{On the restriction $\alpha>\frac{2}{3}$}. This restriction mainly results from the proof technique and is unlikely to be a conceptually fundamental barrier. 
Carefully checking the computations in case (B) of the proof of Theorem \ref{mainstructure}, we see that $\alpha>\frac{2}{3}$ is the threshold value for the whole argument to work. It might be possible to improve the restriction on $d_n$ by a $\sqrt{\log n}$ factor, but getting towards $\alpha<\frac{2}{3}$ may require complete new ideas. Indeed, the inverse Littlewood-Offord theorem on finite fields is well known (see \cite{campos2021singularity}, \cite{campos2022singularity}, \cite{ferber2021counting}) to work well up to $p\leq \exp(\sqrt{n\log n})$ where $n$ is the dimension, which leads to the constraint $\alpha>\frac{2}{3}$ in our band matrix case.

Ideally, we expect our method to work well for all $\alpha\in(\frac{1}{2},1)$, and the range $\alpha\in(\frac{1}{2},\frac{2}{3})$ is only reachable when we can prove a much more refined version of inverse Littlewood-Offord inequality than the one presented here. Invertibility for small bandwidth $\alpha\in(0,\frac{1}{2}]$ is also likely to hold but is clearly not within reach by this method.

\textbf{On the choice of finite field size $p$}. There have been many papers studying integral random matrices over a finite field, see for example \cite{fulman2015stein}, \cite{nguyen2022random} and references therein. This paper appears to be the first to consider the rank of a truly structured integral band matrix over a finite field (but see the recent work of Mészáros \cite{meszaros2024phase} on a random band matrix whose nonzero entries are uniform generated over the $p$-adic integers, and the behavior of the model in \cite{meszaros2024phase} is very different from the model in this paper). In all previous papers considering a square (symmetric) or rectangular matrix on a finite field, the law of the rank of these random matrices often has a universal limiting distribution (see \cite{ferber2023random} \cite{eberhard2022characteristic}) for all prime $p$ growing not too fast in $n$, and in particular we can prove asymptotic non-singularity of these matrices by taking $p\to\infty$ arbitrarily slowly. However, in our band matrix setting, we cannot gain any information on the invertibility of $A_n$ when $p$ does not grow fast enough in $n$, and we also cannot work through the proof when $p$ grows too quickly in $n$. Thus, one really needs a careful choice of $p$ as a function in $n$ in the band matrix setting. 

\end{remark}

Then we begin with the proof of Theorem \ref{mainstructure}.

\begin{proof}[\proofname\ of Theorem \ref{mainstructure}]
    Let $\mu\in(0,\frac{1}{4}]$ be a constant to be determined later, and $v\in\mathbb{F}_p^n\setminus\{0\}$. For each interval $I_k$, we need to study the arithmetic properties of $v_{I_k}$.

    \textbf{Classification of vectors.}
    We will classify the vector $v_{I_k}$ into three possible cases, either (1) that $v_{I_k}$ has a very small support; or (2) $v_{I_k}$ has large support and very nice anti-concentration property, in the sense that $\rho_\mu(v_{I_k})\in(0,\frac{2}{p}]$; or (3) $v_{I_k}$ has a large support but not good enough anti-concentration property, in the sense that $\rho_\mu(v_{I_k})\in[\frac{2}{p},1]$. For the last case, we will consider a dyadic decomposition of the range $\rho_\mu(v_{I_k})\in[\frac{2}{p},1]$. To summarize, we will classify the vectors $v_{I_k}$ into the following three possible categories, for each $k$ with $1\leq k\leq \lfloor n/e_n\rfloor+1$:
    \begin{enumerate}
        \item If $0<|v_{I_k}|\leq K$ where $K$ is given in Lemma \ref{lemmafourier}, then we only need at most $K$ coordinates to define $v_{I_k}$ .
        \item If $|v_{I_k}|\geq K$ and $\rho_\mu(v_{I_k})\leq\frac{2}{p}$, then we can enumerate all vectors $v_{I_k}$ as they can take all the $p^{|I_k|}$ possible vectors.
        \item In the remaining case where $|v_{I_k}|\geq K$ and $\rho_\mu(v_{I_k})\in[\frac{2}{p},1]$, we consider a dyadic decomposition for the interval $[\frac{2}{p},1]\subset\cup_{i=1}^{\log_2 p}[2^{-i},2^{-i+1}]$ for the range of $\rho_\mu(v_{I_k})$. 
    \end{enumerate}
    The total number of intervals needed for the dyadic decomposition is \begin{equation}\label{intervalsdyadicdecom}\#\{\text{Intervals for dyadic decomposition} \}\leq O(\log_2 p)^{n/e_n}=2^{o(n)}\end{equation} under the assumption that $e_n\geq n^\alpha$ and $\alpha>\frac{2}{3}$. 

\textbf{Classification of the arithmetic condition of each vector.}
The proof of this theorem follows from exhausting all the possible values for $\rho_\mu(v_{I_k}),1\leq k\leq \lfloor n/e_n\rfloor+1$ and using the inverse Littlewood-Offord counting theorem, Theorem \ref{inversetheorem}, for vectors in case (3). We will consider two special scenarios (A), (B) and one general scenario (C).

For scenario (A) we assume that all the components $v_{I_k}$ belong to case (3) above except the label $k=n_\mathbf{I}$. For scenario (B) we roughly assume that each $v_{I_k}$ is in either case (1) or case (2) above, and finally the general scenario (C) roughly covers the case where each $v_{I_k}$ can be of category either (1) or (2) or (3). For technical reasons we will consider a slightly more general scenario $(B')$, so that $(A)\cup (B)\cup (B')\cup (C)$ covers all possibilities.

We analyze each scenario in more details below. We also condition on the value $c_{ij}$ in the decomposition $a_{ij}=c_{ij}+u_{ij}$ for all $|i-j|\leq d_n$.

(A) Assume that $\rho_\mu(v_{I_k})\leq\frac{2}{p}$ for all $1\leq k\leq \lfloor n/e_n\rfloor+1:k\neq n_\mathbf{I}$. That is, all components $v_{I_k}$ of $v$, except the $k=n_\mathbf{I}$-th component that we ask about, have very good anti-concentration properties. Then we can upper bound $\mathbb{P}(A_n^\mathbf{I}v=0)$ as follows: for each row of $A_n^\mathbf{I}$ whose index lies in $I_k$, $k\neq n_\mathbf{I}$, we use the randomness in the $D_k$ component in \eqref{blockformofd} and condition on all other random variables in the same row. For the rows with index lying in $I_{n_\mathbf{I}}$, we either use the randomness in $U_{k+1}$ or the randomness in $T_{k-1}$ \footnote{we can choose either of $U_{k+1}$ or $T_{k-1}$ if $k\neq 1$ and $k\neq \lfloor n/e_n\rfloor+1$, but we can only choose $U_2$ for $k=1$ and we can only choose $T_{\lfloor n/e_n\rfloor}$ for $k=\lfloor n/e_n\rfloor+1$}. Then 
$$
\mathbb{P}(A_n^\mathbf{I}v=0)\leq \prod_{k\leq \lfloor n/e_n\rfloor+1,k\neq n_\mathbf{I}}\rho(v_{I_k})^{|I_k|}\cdot \max_{k\in\{n_\mathbf{I}-1,n_\mathbf{I}+1\}}\rho(v_{I_k})^{|I_{n_\mathbf{I}}|-1}\leq  (\frac{2}{p})^{n-1},
$$ where we use $\rho(v_{I_k})\leq\rho_\mu(v_{I_k})$.
Let $\tau>0$ be another constant to be fixed later, and we upper bound the cardinality of all vectors $v_{I_{n_\mathbf{I}}}\in \mathbb{F}_p^{|I_{n_\mathbf{I}}|}$ with $\rho_\mu(v_{I_{n_\mathbf{I}}})\geq\exp(-\tau n^{\alpha/2})$ by 
$$\#\{\text{Candidates for } v_{I_{n_\mathbf{I}}}\}\leq 2^{|I_{n_\mathbf{I}}|}
p^{\frac{2\tau}{\mu} n^{\alpha/2}}\cdot p^{\frac{2\tau}{\mu} n^{\alpha/2}}\cdot (256\exp(\tau n^{\alpha/2}))^{|I_{n_\mathbf{I}}|}
$$ where we apply Theorem \ref{inversetheorem} to get the above upper bound. More precisely, by Theorem \ref{inversetheorem}, we can set up a mapping $f$ from $v_{I_{n_\mathbf{I}}}$ to a subset $T\subset I_{n_\mathbf{I}}$ . The $2^{|I_{n_\mathbf{I}}|}$ factor upper bounds the choice for $T$,  the first $p^{\frac{2\tau}{\mu} n^{\alpha/2}}$ term upper bounds the number of vectors $v_T$, the second term of $p^{\frac{2\tau}{\mu} n^{\alpha/2}}$
upper bounds the number of choice of $v$ in at most $D:=\frac{2\tau}{\mu}n^{\alpha/2}$ coordinates not included in $N_\mu(w)$ with $w:=v_T$, and 
the last term bounds the number of choices of coordinates of $v_i\in N_\mu(w)$. (The case where $|v|\leq D$, which is excluded by Theorem \ref{inversetheorem}, is also upper bounded by the above expression.)

Then combining the previous two expressions and applying a union bound, we get
\begin{equation}\label{case12}\begin{aligned}&\mathbb{P}( \text{ There exists }  v\in\ker A_n^\mathbf{I}: \rho_\mu(v_{I_k})\leq\frac{2}{p} \forall k\neq n_\mathbf{I},\quad  \rho_\mu(v_{I_{n_\mathbf{I}}})\geq\exp(-\tau n^{\alpha/2}))\\&\leq 
2^np^{-n+1}\cdot p^{n-|I_{n_\mathbf{I}}|}\cdot 2^{O_{\mu,\tau}(n^{\alpha/2}\log n)}\exp(\tau n^{\alpha/2}|I_{n_\mathbf{I}}
|),\end{aligned}\end{equation}
where $a=O_{\mu,\tau}(b)$ means $a\leq Cb$ for $C$ a constant depending only on $\mu, \tau$ and $\alpha$. 

Recall that we take $p=\exp(\rho n^{\alpha/2})$ and $|I_{n_\mathbf{I}}|\geq n^\alpha$ with $\alpha>\frac{2}{3}$, we see that as long as $\tau<\rho/2$, the above expression is super-exponentially small, which is actually at most $\exp(-\Omega(n^{1.5\alpha}))$. Thus with probability at least $1-\exp(-\Omega(n^{1.5\alpha}))$, no such vector as described in the first line of \eqref{case12} lies in $\ker A_n^\mathbf{I}$. Formally, we have proven that
\begin{equation}
\label{fa382fa382}    \mathbb{P}(\exists v\in\ker A_n^\mathbf{I}:\rho_\mu(v_{I_k})\leq \frac{2}{p}\forall k\neq n_\mathbf{I},\rho_\mu(v_{I_{n_\mathbf{I}}})\geq\exp(-\tau n^{\alpha/2}))\leq\exp(-\Omega(n^{1.5\alpha})). 
\end{equation}

To summarize, we have proven the following claim:
\begin{Claim}\label{claim01}
      For any $\mu\in(0,\frac{1}{4}]$, any $\rho>0$, and any $\tau<\rho/2$, estimate \eqref{fa382fa382} holds.
\end{Claim}

(B) Then we consider the second case where none of the $v_{I_k}$, $1\leq k\leq \lfloor n/e_n\rfloor+1$ satisfies $\rho_\mu(v_{I_k})\leq\frac{e_n}{p}$. We first assume that $|v_{I_k}|\geq K$ for each $1\leq k\leq \lfloor n/e_n\rfloor+1$  (where $K>0$ is a fixed constant specified in Lemma \ref{lemmafourier}), and treat the other case (that is, where $|v_{I_k}|\leq K$ for some $k$), at the end of scenario (B) by a simple modification. Let $t_1,\cdots,t_{\lfloor n/e_n\rfloor+1}$ be integers such that $\rho_\mu(v_{I_k})\in[2^{-t_k},2^{-t_k+1}]$ for each $k$. Then we upper bound the probability $\mathbb{P}(A_n^\mathbf{I}v=0)$ by conditioning on all the random variables away from the diagonal blocks $D_1,D_2,\cdots,D_{\lfloor n/e_n\rfloor+1}$ and using randomness only from these diagonal blocks:
$$\begin{aligned}&
\mathbb{P}(A_n^\mathbf{I}v=0)\leq\prod_{k=1}^{\lfloor n/e_n\rfloor+1}|\rho(v_{I_k})|^{|I_k|-1_{k=n_\mathbf{I}}}\\&\leq c^{n-1}\prod_{k=1}^{\lfloor n/e_n\rfloor+1}|\rho_\mu(v_{I_k})|^{|I_k|-1_{k=n_\mathbf{I}}}\leq \frac{p}{c}\cdot c^n\prod_{k=1}^{\lfloor n/e_n\rfloor+1}(2^{-t_k+1})^{|I_k|},
\end{aligned}$$ where in the second line we take any value of $c\in(0,1)$ to be fixed at the end of proof and use Lemma \ref{lemmafourier} to get a value $\mu\in(0,\frac{1}{4}]$ such that Lemma \ref{lemmafourier} holds with this value $c$.  The indicator function $1_{k=n_\mathbf{I}}$ comes from removal of the $\mathbf{I}$-th row from $A_n^\mathbf{I}$.

Then we compute the cardinality of vectors $v_{I_1},\cdots,v_{I_k},\cdots,v_{I_{\lfloor n/e_n\rfloor+1}}$ satisfying the stated constraint on anticoncentration. Set $D_k=\frac{2}{\mu}(t_k+1)$ for each $k$, as we apply Theorem \ref{inversetheorem}, for each $v_{I_k}$ we have at most $2^{|I_k|}$ choices for $T$, $p^{D_k}$ choices for $w$, and $2^{(t_k+9)}$ choices for every other coordinate of $v_{I_k}$. By Theorem \ref{inversetheorem}, the total cardinality for these vectors of $v_{I_k}$ is at most 
\begin{equation}\label{balloonsinthesky}\#\{\text{Candidates for } v_{I_k}\}\leq 
2^{|I_k|+9}p^{2D_k}(2^{\mu D_k/2})^{|I_k|-D_k}\leq 2^{|I_k|+9}\cdot2^{(t_k+1)|I_k|}\cdot e^{\frac{8}{\mu}(\log p)^2},
\end{equation} where the first term $2^{|{I_k}|}$ comes from the choice of $T\subset I_k$, the two powers $p^{D_k}$ come from the choice of $v_{ T}$ and the coordinates of $v_{I_k}$ not in $T$, and the last term is the cardinality of $N_\mu(v_T)$. Finally, note that $t_k+1\leq\log_2 p\leq 2\log p$, which leads to the $\frac{8}{\mu}(\log p)^2$ exponent.

Taking the product for all $k=1,2,\cdots,\lfloor n/e_n\rfloor+1$, we have via a direct union bound,
\begin{equation}\label{thesecondline}\begin{aligned}&
\mathbb{P}(\text{There exists }v\in \ker A_n^\mathbf{I}:\rho_\mu(v_{I_k})\in[2^{-t_k},2^{-t_k+1}]\forall k,|v_{I_k}|\geq K)\\&\leq pc^{n-1}2^{20n}e^{\frac{8}{\mu}(\log p)^2\frac{n}{e_n}}.\end{aligned}
\end{equation} By our assumption, $(\log p)^2\frac{n}{e_n}\leq \rho^2n^\alpha\frac{n}{n^\alpha/s}\leq\rho^2 sn$ for $s$ fixed (where $e_n=\lfloor d_n/s\rfloor$), so that by setting $c>0$ small [this fixes $\mu$] and then setting $\rho>0$ sufficiently small relative to $\mu$, we can guarantee that the second line of \eqref{thesecondline} is at most $\exp(-\Omega(n))$. Since the number of dyadic decompositions is at most $2^{o(n)}$ by \eqref{intervalsdyadicdecom}, the eventual probability for a vector $v$ in case (B) to be in $\ker A_n^\mathbf{I}$ (such that $|v_{I_k}|\geq K$ for each $k$) is $\exp(-\Omega(n))$.
That is, 
\begin{equation}\label{equations366unmatch}
    \mathbb{P}(\exists v\in\ker A_n^\mathbf{I}\setminus\{0\}:\rho_\mu(v_{I_k})\geq \frac{e_n}{p}\forall k,|v_{I_k}|\geq K\forall k)\leq\exp(-\Omega(n)). 
\end{equation}

 We now treat the leftover case where for some $k$, we have $v_{I_k}=0$ or $0<|v_{I_k}|<K$ for the constant $K>0$ in Lemma \ref{lemmafourier}. Denote by $\mathcal{R}$ the subset of indices $k$ such that $|v_{I_k}|\leq K$. Then 
 whenever $0<|v_{I_k}|\leq K$ we use the trivial fact that $\rho(v_{I_k})\leq\frac{1}{2}$ combined with the following upper bound on the number of such $v_{I_k}$: the cardinality is bounded by $$\#\{\text{Candidates for } v_{I_k}:|v_{I_k}|\leq K\}\leq(e_n)^Kp^K=e^{O_K(n^{\alpha/2}\log n)}.$$ If $v_{I_k}=0$, then we proceed trivially. The similarly as in estimate \eqref{equations366unmatch}, we deduce that for any possible choice of the set $\mathcal{R}$, we have
 \begin{equation}\label{A61345}
    \mathbb{P}(\exists v\in \ker A_n^\mathbf{I}:\rho_\mu(v_{I_k})\geq\frac{e_n}{p}\forall k,|v_{I_k}|\geq K\text{ if and only if } k\notin\mathcal{R})\leq\exp(-\Omega(n)).
\end{equation}
Enumerating all the possible subsets $\mathcal{R}\subset[\lfloor n/e_n\rfloor+1]$, we have considered all vectors in our type (B), and thus we have verified that \begin{equation}
    \mathbb{P}(\exists v\in\ker A_n^\mathbf{I}\setminus\{0\}:\rho_\mu(v_{I_k})\geq \frac{e_n}{p}\forall k)\leq\exp(-\Omega(n). 
\end{equation}

(B') We also need to consider the case where for exactly one $k$ we have $\rho_\mu(v_{I_k})\leq \frac{e_n}{p}$ holds. For this case we use the same estimate as in case (B) for $v_{I_l},l\neq k$, but for component $v_{I_k}$ we use all the $p^{|I_k|}$ possible values of $v_{I_k}$. Then the probability estimate in \eqref{thesecondline} would now have the form $pc^{n-1-|I_k|}2^{20n}e^{\frac{8}{\mu}(\log p)^2\frac{n}{e_n}}\cdot (e_n)^{e_n}$. But we have $(e_n)^{e_n}=2^{o(n)}$ (this follows from our assumption $e_n\leq n^{1-\beta}$), so the above quantity is again $\exp(-\Omega(n))$. That is, we can verify that 
\begin{equation}
\label{before2.6}    \mathbb{P}(\exists v\in\ker A_n^\mathbf{I}\setminus\{0\}:\rho_\mu(v_{I_k})\geq \frac{e_n}{p}\text{ for all but one } k\in[\lfloor n/e_n\rfloor+1])\leq\exp(-\Omega(n). 
\end{equation}

To summarize, we have proven the following claim:
\begin{Claim}\label{claim02}
  We can find a constant $c_0>0$ such that for any $c\in(0,c_0)$, we can find $\mu=\mu(c)\in(0,\frac{1}{4}]$ such that, whenever $\rho\leq\rho(\mu)$ [where $\rho(\mu)$ is a fixed function of $\mu$ only], then estimate \eqref{before2.6} holds.
\end{Claim}

(C) Finally, we consider the general case where $\rho_\mu(v_{I_{n_\mathbf{I}}})\geq\exp(-\tau n^{\alpha/2})$, with $\tau=\rho/2$, and there exists some $k\neq n_{\mathbf{I}},k\neq \lfloor n/e_n\rfloor+1$ such that $\rho_\mu(v_{I_k})\leq\frac{e_n}{p}.$ That is,
\begin{equation}
\label{436expeditions}    \exists k\neq n_\mathbf{I}:\rho_\mu(v_{I_k})\leq\frac{e_n}{p},\quad\rho_\mu(v_{I_{n_\mathbf{I}}})\geq \exp(-\tau n^{\alpha/2}).
\end{equation}
Note that $\rho_\mu(v_{I_k})\leq\frac{e_n}{p}$ immediately implies $ |v_{I_k}|\geq K$ for any fixed $K$ when $n$ is large enough. Observe that the complement of case (C) is studied in case (B) and (B'). Fix one such value $k_*$ such that $\rho_\mu(v_{I_{k_*}})\leq\frac{e_n}{p}$. We assume without loss of generality that $k_*<n_\mathbf{I}$.

Then we can find a sequence of integers $n_\mathbf{I}=k_0>k_1>k_2\geq\cdots>k_\Delta:=k_*$ ($\Delta\in\mathbb{N}$ is an integer) such that $k_i$ is the largest index in $[1,k_{i-1}]$ such that $\rho_\mu(v_{I_{k_i}})\leq\rho_\mu(v_{I_{k_{i-1}}})$ for each $i$, where we denote by $k_0=n_\mathbf{I}$. That is, we find the maximal chain of indices $k_0,k_1,\cdots,k_\Delta$ where $\rho_\mu(v_{I_{k_i}})$ is decreasing as $i$ increases from $0$ to $\Delta$. We clearly have by definition of $n_\mathbf{I}=k_0$ and $k_*=k_\Delta$ that
\begin{equation}\label{quotients1}\exp(-\rho n^{\alpha/2}/4)\geq
\frac{\rho_\mu(v_{I_{k_*}})}{\rho_\mu(v_{I_{n_\mathbf{I}}})}=\prod_{i=1}^{\Delta}\frac{\rho_\mu(v_{I_{k_i}})}{\rho_\mu(v_{I_{k_{i-1}}})}\geq\prod_{i=1}^{\Delta}\frac{\rho_\mu(v_{I_{k_i}})}{\rho_\mu(v_{I_{k_{i}+1}})}.
\end{equation}
We upper bound $\mathbb{P}(A_n^\mathbf{I}v=0)$ via the following approach: for rows indexed in $I_{k_i+1}$ we use randomness from $T_{I_{k_i}}$ (see \eqref{blockformofd}), and for all other rows indexed by $I_l$ we use randomness from the diagonal component $D_l$. This helps us to make the best use of anticoncentration properties for each component of $v$. More formally, for this vector $v$,
$$
\mathbb{P}(A_n^\mathbf{I}v=0)\leq \prod_{i\in\lfloor n/e_n\rfloor+1:i\notin k_1+1,\cdots, k_\Delta+1\}} \sup_w\mathbb{P}(D_iv_{I_i}=w)
\cdot\prod_{i=1}^\Delta \sup_w\mathbb{P}(T_{I_{k_i}}v_{I_{k_i}}=w).$$
Then we have, using Fact \ref{fact2.111}\footnote{To clarify, in contrast to case (B), in case (C) we only use the trivial monotonicity inequality $\rho(v_{I_k})\leq\rho_\mu(v_{I_k})$ instead of the more refined one $\rho(v_{I_k})\leq c\rho_\mu(v_{I_k})$, and the former always holds for any $v_{I_k}$.} which states $\rho(v_{I_k})\leq\rho_\mu(v_{I_k})$ for any $k$ and $\mu\in(0,\frac{1}{4}]$,
\begin{equation}\label{quotients2}
\mathbb{P}(A_n^\mathbf{I}v=0)\leq  
\prod_{i=1}^{\lfloor n/e_n\rfloor+1}|\rho_\mu(v_{I_i})|^{|I_i|}\cdot \prod_{i=1}^\Delta \left|\frac{\rho_\mu(v_{I_{k_i}})}{\rho_\mu(v_{I_{k_{i}+1}})}\right|^{|I_{k_i+1}|}.\end{equation}

We now count the number of possible vectors $v$ that contribute to $A_n^\mathbf{I}v=0$. Let $\Omega_1,\Omega_2$ be a partition of $\lfloor n/e_n\rfloor+1$ such that $k\in\Omega_1$ if $\rho_\mu(v_{I_k})\geq\frac{2}{p}$, and $k\in\Omega_2$ if $\rho_\mu(v_{I_{k}})\leq\frac{2}{p}$. Then for each $k\in\Omega_1$ we can find $t_k\in\mathbb{N}$ such that $\rho_\mu(v_{I_k})\in [2^{-t_k},2^{-t_k+1}]$. Then applying Theorem \ref{inversetheorem}, we get similarly as in \eqref{balloonsinthesky} that
$$
k\in\Omega_1,\rho_\mu(v_{I_k})\in [2^{-t_k},2^{-t_k+1}]\Rightarrow \#\{\text{Candidates for } v_{I_k}\}\leq 2^{|I_k|+9}\cdot2^{(t_k+1)|I_k|}\cdot e^{\frac{8}{\mu}(\log p)^2}.
$$ If $k\in\Omega_2$ then we simply bound the cardinality of possible vectors $v_{I_k}$ by $p^{|I_k|}$.

Then by a direct union bound, we have for each choice of $\{t_i\}$ and $\Omega_1,\Omega_2$,
$$\begin{aligned}&
\mathbb{P}(\text{There exists }v\in\ker A_n^\mathbf{I}:\rho_\mu(v_{I_k})\in[2^{-t_k},2^{-t_k+1}]\forall k\in\Omega_1,\rho_\mu(v_{I_k})\leq\frac{2}{p}\forall k\in\Omega_2,\text{  }\eqref{436expeditions})\\&\leq 2^{20n}e^{\frac{8}{\mu}(\log p)^2\frac{n}{e_n}}\prod_{k:\rho_\mu(v_{I_k})\geq\frac{2}{p}}2^{(t_k+1)|I_k|}\prod_{k:\rho_\mu(v_{I_k})\leq\frac{2}{p}}p^{|I_k|}\\&\quad \cdot \prod_{k=1}^{\lfloor n/e_n\rfloor+1}|\rho_\mu(v_{I_k})|^{|I_k|}\cdot \exp(-\rho n^{\alpha/2}|I_{n_\mathbf{I}}|/4),
\end{aligned}$$

where we have applied \eqref{quotients1} to the last term on the right hand side of \eqref{quotients2}. Simplifying the above expression and noting that $n^{\alpha/2}|I_{n_\mathbf{I}}|\geq n^{1.5\alpha}$ and $\alpha>\frac{2}{3}$, the last display is at most $\exp(-\Omega(n^{1.5\alpha}))$. Summing over all possible choices of $\Omega_1,\Omega_2$ and $\{t_i\}$, which has cardinality $2^{o(n)}$, we conclude that
\begin{equation}
 \label{estimate46789}   \mathbb{P}(\text{There exists }v\in\ker A_n^\mathbf{I}\text{ satisfying }\eqref{436expeditions}\text{ for some }n_\mathbf{I},k_*)\leq\exp(-\Omega(n^{1.5\alpha})).
\end{equation}
To summarize, we have proven the following claim:
\begin{Claim}\label{claim03}
    For any $\mu\in(0,\frac{1}{4}]$ and $\rho>0$,  taking $\tau=\rho/2$, then estimate \eqref{estimate46789} holds.
\end{Claim}

Summarizing cases (A), (B), (B') and (C), which result in Claim \ref{claim01}, \ref{claim02} and \ref{claim03}, we conclude that
there exists a suitable choice of $\rho,\tau>0$ and $\mu\in(0,\frac{1}{4}]$ such that
\begin{equation}
 \mathbb{P}(\text{There exists }v\in\ker A_n^\mathbf{I}:\rho_\mu(v_{I_{n_\mathbf{I}}})\geq\exp(-\tau n^{\alpha/2})\text{ for some }n_\mathbf{I})   \leq \exp(-\Omega(n))
.\end{equation}
 This completes the entire proof of Theorem \ref{mainstructure}.\end{proof}

\section{Invertibility for the symmetric band matrix model}

This section is devoted to the proof of Theorem \ref{symmetricmaintheorem}. We divide the proof into two separate parts: a structural characterization for null vectors of principal submatrices in Section \ref{sections3.1}, and then the linear algebraic reductions which complete the proof of Theorem \ref{symmetricmaintheorem} in Section \ref{section3.223}.

\subsection{Justifying shapes of null vectors of submatrices}\label{sections3.1}

We first recall the following constructions in Section \ref{section2.2}. We define $e_{n}=\lfloor d_{n}/s\rfloor$
(for some $s\in(3,6)$ to be fixed) and consider a partition of $[1,n]$ by intervals of length $e_n$:
\begin{equation}\label{911i1i2}
[1,n]=\cup_{k=1}^{\lfloor n/e_{n}\rfloor +1}I_k,\quad I_k=\begin{cases} [(k-1)e_n+1,ke_n],\quad \forall 1\leq k\leq \lfloor n/e_n\rfloor \\ [\lfloor n/e_n\rfloor e_n+1,n],\quad k=\lfloor n/e_n\rfloor+1,
\end{cases}\end{equation} where we choose some $s\in(3,6)$ such that $n-e_n\lfloor n/e_n\rfloor\geq e_n/2$. 

\begin{theorem}\label{theoremsymmetriccase}
    Let $S_n$ be the random matrix defined in Theorem \ref{symmetricmaintheorem}. Then we can find a choice of $\rho>0$ and $\mu\in(0,\frac{1}{4}]$ such that, with $p=\exp(\rho n^{\alpha/2})$ and let $\ker_{\mathbb{F}_p} S_n$ denote the kernel of $S_n$ modulo $\mathbb{F}_p$, 
    $$\begin{aligned}&\mathbb{P}(\forall v\in \ker_{\mathbb{F}_p} S_n\setminus\{0\},\forall 1
    \leq k\leq \lfloor n/e_n\rfloor:\min(\rho_\mu(v_{I_k}),\rho_\mu(v_{I_{k+1}}))\leq\exp(-\tau n^{\alpha/2}))\\&\geq 1-\exp(-\Omega(n)).
    \end{aligned}$$
\end{theorem}

In plain words, Theorem \ref{theoremsymmetriccase} proves that, with very high probability, all the nonzero kernel vectors to $S_n$ on $\mathbb{F}_p$ must satisfy that, at least one of its \textbf{any} two consecutive components $v_{I_k},v_{I_{k+1}}$ must have good anti-concentration properties in terms of $\rho_\mu(\cdot)$.
    
We make some preparations for the proof of Theorem \ref{theoremsymmetriccase}. We define the following matrix from $S_n=(s_{ij})$ that will be used for comparison: 
\begin{equation}\label{symremovedsubmatrix}\begin{aligned}
\mathcal{D}&:=\{s_{ij}\mathbf{1}_{i,j\in I_k\text{ or }i\in I_{k-1},j\in I_k\text{ or } i\in I_k,j\in I_{k-1}\text{ for some }k=1,\cdots,\lfloor n/e_n\rfloor+1}\}.\end{aligned}
\end{equation} (where we interpret $I_0=\emptyset$).

In matrix form, we can write $\mathcal{D}$ in the following block form:
\begin{equation}\label{symblockformofd}
    \mathcal{D}=\begin{bmatrix}
        D_1&U_2&0&\cdots&U_1^T\\U_2^T&D_2&U_3&\cdots&0\\0&U_3^T&D_3&\cdots&0\\\ \vdots&\vdots&\vdots&\ddots&U_{\lfloor n/e_n\rfloor+1}\\U_1&0&0&U_{\lfloor n/e_n\rfloor+1}^T&D_{\lfloor n/e_n\rfloor+1}
    \end{bmatrix}.
\end{equation} By the assumption of Theorem \ref{symmetricmaintheorem}, the entries in the blocks $U_i,1\leq i\leq \lfloor n/e_n+1\rfloor$, are independent and uniformly distributed over $\{\pm 1\}+c_{ij}$ for some fixed location-dependent constant $c_{ij}\in\mathbb{Z}$ that do not matter to us. Meanwhile, the diagonal blocks $D_i,1\leq i\leq \lfloor n/e_n\rfloor+1$ are \textbf{symmetric} random matrices with entries distributed on $\{\pm 1\}+c_{ij}$ for some fixed location-dependent constant $c_{ij}$. In other words, $\mathcal{D}$ is a matrix extracted from $S_n$ by zeroing out all the entries not in the main diagonal and first off-diagonal blocks labeled by $I_k,k=1,\cdots,\lfloor n/e_n\rfloor+1$. Note that $U_1\neq 0$, as we assume our variance profile wraps around, see Remark \ref{cyclicwarparound}. Thus the variance profile of $\mathcal{D}$ is translation invariant.

To further exploit independence, we need to use the following structural decomposition of diagonal blocks $D_i$: for each $1\leq i\leq \lfloor n/e_n\rfloor+1$,
\begin{equation}\label{decompositionoftheblockd}
D_i=\begin{bmatrix}
    (D_i)_{S\times S}&(D_i)_{S\times S^c}\\(D_i)_{S^c\times S}&(D_i)_{S^c\times S^c}
\end{bmatrix},
\end{equation} where $S$ is any subset of the rows/columns of $D_i$, and $S^c=I_i\setminus S$.

The proof of Theorem \ref{theoremsymmetriccase} is based on analyzing all possible types of vectors in $\ker S_n$. We will categorize all kernel vectors into four structure types: type I, II, III and IV. We will analyze each structure separately in each subsection, after which the proof of Theorem \ref{theoremsymmetriccase} directly follows.

\subsubsection{Eliminating kernel vectors of structure type I}

This part corresponds to scenario (A) in the proof of Theorem \ref{mainstructure}. 

\begin{notation} Fix any constant $\mu\in(0,1]$.
    We say that a vector $v\in\mathbb{F}_p^n$ is of type I$(\mu)$ if,
for some index $n_\mathbf{I}$ with $1\leq n_\mathbf{I}\leq \lfloor n/e_n\rfloor+1$ we have $\rho_\mu(v_{I_k})\leq\frac{2}{p}$ for all $k\neq n_\mathbf{I}$ and $\rho_\mu(v_{I_{n_\mathbf{I}}})\geq\exp(-\tau n^{\alpha/2})$. Here $p=\exp(\rho n^{\alpha/2})$, and $\rho,\tau$ are constants that will be fixed later. 
\end{notation}

\begin{remark}
    The definition of type I$(\mu)$ depends on the parameter $\mu>0$. When the context is clear, we will omit the $(\mu)$ symbol and simply say $v$ is of type I.
\end{remark}

\begin{lemma}\label{lemma3.9}
   We can find a constant $c_0>0$ such that for any $c\in(0,c_0)$, we can find $\mu=\mu(c)\in(0,\frac{1}{4})$ such that, whenever $\rho\leq\rho(\mu)$ [where $\rho(\mu)$ is a fixed function of $\mu$ only] and $\tau=\rho/32$, the following holds with probability $1-\exp(-\Omega(n))$: There exists no vector $v\in\mathbb{F}_p^n\setminus\{0\}$ of type I$(\mu)$, such that $S_n\cdot v=0$. 
 \end{lemma}
 \begin{remark}
     Actually, in the proof of Lemma \ref{lemma3.9} , we can take any $c_0,$ any $\mu\in(0,\frac{1}{4})$ and any $\rho>0$. We make these additional constraints on $c,\rho,\mu$ as they will show up when we consider vectors of type $II,III,IV$ later.

 \end{remark}

 \begin{proof}

     Without loss of generality, we can assume $n_\mathbf{I}=1$ and the other cases are similar. Then for any nonzero vector $v$ of type I, we use the following method to upper bound $\mathbb{P}(S_n\cdot v=0)$. It suffices to bound $\sup_{w\in\mathbb{F}_p^n}\mathbb{P}(\mathcal{D}\cdot v=w)$, then for rows labeled by $I_1$ we use randomness in $U_2$, for rows labeled by $I_2$ we use randomness in $U_3$, all the way down to the rows associated to $I_{\lfloor n/e_n\rfloor}$ we use randomness in $U_{[\lfloor n/e_n\rfloor+1]}$. Finally, for the last block row associated to $I_{\lfloor n/e_n\rfloor+1]}$ we use the randomness in $D_{\lfloor n/e_n\rfloor}+1$. More precisely, we will use an iterative conditioning argument that is written in formula as follows:

     Let  $(D)_1^{k-1}$ denote the subset $\{D_1,\cdots,D_{k-1}\}$ and $(U)_1^k$ denote the subset $\{U_1,\cdots,U_k\}$, then  $$\begin{aligned}&
\mathbb{P}(S_n\cdot v=0)\leq\sup_w\mathbb{P}(\mathcal{D}\cdot v=w)
\\&\leq\prod_{k=1}^{\lfloor n/e_n\rfloor+1}\sup_{w,(D)_1^{k-1},(U)_1^k}\mathbb{P}(U_{k}^Tv_{I_{k-1}}+D_{k}v_{I_k}+U_{{k+1}}v_{I_{k+1}}=w\mid (D)_1^{k-1},(U)_1^k )
\\&\leq \left(\prod_{k=2}^{\lfloor n/e_n\rfloor+1}\sup_w \mathbb{P}(U_k v_{I_k}=w)\right)\cdot\sup_w\mathbb{P}(D_{\lfloor n/e_n\rfloor+1}\cdot v_{I_{\lfloor n/e_n\rfloor+1}}=w),
    \end{aligned} $$
where we write $I_{\lfloor n/e_n\rfloor+2}=I_1$.
Since each $U_k$ has independent entries, we can bound for each $2\leq k\leq \lfloor n/e_n\rfloor+1$, $$\sup_w\mathbb{P}(U_kv_{I_k}=w)\leq \rho(v_{I_k})^{|I_{k-1}|}\leq \rho_\mu(v_{I_k})^{|I_{k-1}|}\leq(\frac{2}{p})^{|I_{k-1}|}$$ by the assumption that the vector $v$ is of type I, and that $\mu\in(0,\frac{1}{4}]$.
For the last term, we decompose $D_{\lfloor n/e_n\rfloor+1}$ as in \eqref{decompositionoftheblockd} by taking $S$ to contain the  smaller half of the integers in $I_{n/e_n+1}$. Under the assumption that $\rho_\mu(v_{I_{\lfloor n/e_n\rfloor}+1})\leq \frac{2}{p}$, let $v^1,v^2$ denote the restriction of $v_{I_{\lfloor n/e_n\rfloor+1}}$ onto $S$ and $S^c$ respectively, we claim that $\min(\rho_\mu(v^1),\rho_\mu(v^2))\leq\sqrt{\frac{2}{p}}$. This follows from the fact that $\rho_\mu(v_{I_{\lfloor n/e_n\rfloor+1}})\geq \rho_\mu(v^1)\rho_\mu(v^2)$, which follows from independence in the definition of $\rho_\mu$. Suppose without loss of generality that $\rho_\mu(v^1)\leq\sqrt{\frac{2}{p}}$, then we can bound 
\begin{equation}\label{howdowecanboundthis}
\sup_w\mathbb{P}(D_{\lfloor n/e_n\rfloor+1}\cdot v_{I_{\lfloor n/e_n\rfloor+1}}=w)\leq\sup_w\mathbb{P}
((D_{\lfloor n/e_n\rfloor+1})_{S^c\times S}\cdot v^1=w)\leq (\sqrt{\frac{2}{p}})^{|I_{\lfloor n/e_n\rfloor+1}|/2}.\end{equation}
 In summary, for this vector $v$, we have
 $$
\mathbb{P}(S_n\cdot v=0)\leq(\frac{2}{p})^{n-|I_1|/2}(\sqrt{\frac{2}{p}})^{|I_1|/4},
 $$where we use the assumption that $|I_{\lfloor n/e_n\rfloor+1}|\geq|I_1|/2$.

 Now we enumerate all the possible vectors $v$ of type I. The choice of $v_{I_k},2\leq k\leq \lfloor n/e_n\rfloor+1$, is arbitrary, so we have $p^{n-|I_1|}$ choices.
    Then we upper bound the cardinality of all vectors $v_{I_1}\in \mathbb{F}_p^{|I_1|}$ with $\rho_\mu(v_{I_1})\geq\exp(-\tau n^{\alpha/2})$ by 
\begin{equation}\label{cardinalityforv1}\#\{\text{Candidates for } v_{I_1}\}\leq 2^{|I_1|}
p^{\frac{2\tau}{\mu} n^{\alpha/2}}\cdot p^{\frac{2\tau}{\mu} n^{\alpha/2}}\cdot (256\exp(\tau n^{\alpha/2}))^{|I_{1}|}
,\end{equation}where this upper bound comes from applying Theorem \ref{inversetheorem}. More precisely, Theorem \ref{inversetheorem} yields a mapping $f$ from $v_{I_I}$ to a subset $T\subset I_1$, $|T|\leq\frac{2}{\mu}\log\rho_\mu(v_{I_1})^{-1}$ and a subset $N_\mu((v_{I_1})_T)\subset\mathbb{F}_p$ of cardinality at most $\frac{256}{\rho_\mu(v)}$ such that $(v_{I_1})_i\in N_\mu((v_{I_1})_T)$ for all but at most $\frac{2}{\mu}\log\rho_\mu(v)^{-1}$ indices. In the
above expression, the $2^{|I_{1}|}$ factor upper bounds the choice for $T$,  the first $p^{\frac{2\tau}{\mu} n^{\alpha/2}}$ term upper bounds the number of vectors $(v_{I_1})_T$, the second such term
upper bounds the number of choice of $v_{I_1}$ in at most $D:=\frac{2\tau}{\mu}n^{\alpha/2}$ coordinates of $I_1$ not included in $N_\mu((v_{I_1})_T)$, and 
the last term bounds the number of choices of coordinates of $v_i\in N_\mu((v_{I_1})_T)$. (The case where $|v|\leq D$, which is excluded by Theorem \ref{inversetheorem}, is also upper bounded by the above expression).

Then combining the previous two expressions and applying a union bound, we get
\begin{equation}\label{case1}\begin{aligned}&\mathbb{P}( \exists   v\in\ker S_n\setminus\{0\}\text{ of type I})\leq 
2^n(\sqrt{\frac{2}{p}})^{|I_1|/4} 2^{O_{\mu,\tau}(n^{\alpha/2}\log n)}\exp(\tau n^{\alpha/2}|I_{1}
|).\end{aligned}\end{equation}
where $a=O_{\mu,\tau}(b)$ means $a\leq Cb$ for $C$ a constant depending only on $\mu, \tau$ and $\alpha$. 

Recall that we take $p=\exp(\rho n^{\alpha/2})$ and $|I_{1}|\geq n^\alpha$ with $\alpha=\frac{2}{3}$. Then we see that as long as $\tau=\rho/32$, the above expression is super-exponentially small, which is actually at most $\exp(-\Omega(n^{1.5\alpha}))$. Thus with probability at least $1-\exp(-\Omega(n^{1.5\alpha}))$, no such vector of type I lies in $\ker S_n$.
\end{proof}

\subsubsection{Eliminating kernel vectors of structure type II}

This part corresponds to scenario (B) in the proof of Theorem \ref{mainstructure}. 

\begin{notation}
    We say that a vector $v\in\mathbb{F}_p^n\setminus\{0\}$ is of type II$(\mu)$ if,
for each $1\leq k\leq \lfloor n/e_n\rfloor+1$, we have $\rho_\mu(v_{I_k})\geq\frac{e_n}{p}$.
\end{notation}

We shall recall the following lemma:

\begin{lemma}\label{lemma3.97}[\cite{campos2022singularity}, Lemma 4] Fix $\mu\in(0,\frac{1}{4}]$ and $n\in\mathbb{N}$. For any $v\in \mathbb{F}_p^n$ we can find $T\subseteq [n]$ such that $v_i\in N(v_T)$ for all $i\notin T$, that $\rho_\mu(v_T)\leq (1-\mu/2)^{|T|}$, and $$|T|\leq\frac{2}{\mu}\log\frac{1}{\rho_\mu(v)}.$$
\end{lemma}

\begin{lemma}\label{newtype2} We can find a constant $c_0>0$ such that for any $c\in(0,c_0)$, we can find $\mu=\mu(c)\in(0,\frac{1}{4}]$ such that, whenever $\rho\leq\rho(\mu)$ [where $\rho(\mu)$ is a fixed function of $\mu$ only],
     the following holds with probability $1-\exp(-\Omega(n))$: There exists no vector $v\in\mathbb{F}_p^n\setminus\{0\}$ of type II$(\mu)$, such that $S_n\cdot v=0$. 
 \end{lemma}
 \begin{proof}\textbf{Step 1: Global reduction.}
We first assume that $|v_{I_k}|\geq K$ for each $1\leq k\leq \lfloor n/e_n\rfloor+1$  (where $K>0$ is a fixed constant specified in Lemma \ref{lemmafourier}). The other case where $|v_{I_k}|\leq K$ for some $k$ is treated via a straightforward modification at the end of proof. 

Since the variance profile of $\mathcal{D}$ is translation invariant, we can assume without loss of generality that $\rho_\mu(v_{I_1})\geq\rho_\mu(v_{I_k})$ for each $k$. Otherwise, we simply permute the rows and columns of $\mathcal{D}$.

Let $t_1,\cdots,t_{\lfloor n/e_n\rfloor+1}$ be integers such that $\rho_\mu(v_{I_k})\in[2^{-t_k},2^{-t_k+1}]$ for each $k$. Then we upper bound the probability $\mathbb{P}(S_n\cdot v=0)$ as in the previous proof:

\begin{equation}\label{caseB1}\begin{aligned}&
\mathbb{P}(S_n\cdot v=0)\leq\sup_w\mathbb{P}(\mathcal{D}\cdot v=w)\\&\leq \left(\prod_{k=2}^{\lfloor n/e_n\rfloor+1}\sup_w \mathbb{P}(U_k v_{I_k}=w)\right)\cdot\sup_w\mathbb{P}(D_{\lfloor n/e_n\rfloor+1}\cdot v_{I_{\lfloor n/e_n\rfloor+1}}=w)\\&\leq \prod_{k=2}^{\lfloor n/e_n\rfloor+1}\rho(v_{I_k})^{|I_{k-1}|} \cdot \sup_w\mathbb{P}(D_{\lfloor n/e_n\rfloor+1}\cdot v_{I_{\lfloor n/e_n\rfloor+1}}=w)
 \\&\leq c^{n-|I_{\lfloor n/e_n\rfloor+1}|}\cdot \prod_{k=2}^{\lfloor n/e_n\rfloor+1}\rho_\mu(v_{I_k})^{|I_{k-1}|} \cdot \sup_w\mathbb{P}(D_{\lfloor n/e_n\rfloor+1}\cdot v_{I_{\lfloor n/e_n\rfloor+1}}=w)   \\&\leq c^{n-|I_{\lfloor n/e_n\rfloor+1}|}\cdot \prod_{k=1}^{\lfloor n/e_n\rfloor}\rho_\mu(v_{I_k})^{|I_{k}|} \cdot \sup_w\mathbb{P}(D_{\lfloor n/e_n\rfloor+1}\cdot v_{I_{\lfloor n/e_n\rfloor+1}}=w)   ,\end{aligned} \end{equation} where in the fourth line, we take any value of $c\in(0,1)$ to be fixed at the end of proof and use Lemma \ref{lemmafourier} to get a value $\mu\in(0,\frac{1}{4}]$ such that Lemma \ref{lemmafourier} holds with this value $c$. The inequality in the fifth line follows from the assumption that $\rho_\mu(v_{I_1})\geq\rho_\mu(v_{I_k})$ for each $k$, and that $|I_k|$ is a constant for all $1\leq k\leq \lfloor n/e_n\rfloor$.

\textbf{Step 2: Analyzing one symmetric matrix block.} The proof in this step is an adaptation of the journal version of \cite{campos2022singularity}.
To simplify the notation, we denote $n_*=\lfloor n/e_n\rfloor+1$ and $I_*=I_{n_*}$. We now upper bound $\sup_w\mathbb{P}(D_{I_*}\cdot v_{I_*}=w)$.
     For each vector $\hat{v}:=v_{I_*}$, by Theorem \ref{inversetheorem} we can find a subset $T\subset I_*$
such that $\hat{v}_i\in N_\nu(\hat{v}_T)$ for all $i\notin T$. That is, let $\mathcal{V}=\{\hat{v}\in\mathbb{F}_p^n\setminus\{0\}:\rho(\hat{v})\geq\beta\}$ where $\beta=\Theta(e_n/p)$, we let $f(\hat{v}):=(T,\hat{v}_T)$ where $T$ is the set obtained by Theorem \ref{inversetheorem}. Let $\mathcal{S}:=f(\mathcal{V})$. Then we can bound 
$$
|\mathcal{S}|\leq 2^n\cdot p^{\frac{2}{\mu}\log p}\leq 2^{(1+\frac{2\rho}{\mu})n^\alpha}.
$$ For each $s=(T,u)\in\mathcal{S}$, we have that for each $\hat{v}\in f^{-1}(s)$ we must have that $\hat{v}_i\in N_\mu(u)$ for all $i\notin T$, and thus by \eqref{doublecounting}, 
$$
|f^{-1}(s)|\leq |N_\mu(u)|^{|I_*|-|T|}\leq (\frac{2}{\rho_\mu(u)})^{|I_*|-|T|}.
$$Then for each 
$\hat{v}\in f^{-1}(s)$, let $T\subset S$ where $|S|=\min\{|T|+K,|\hat{v}|\}$ and that $\hat{v}_i\neq 0$ for all $i\in S$. Then we can bound 
$$
\mathbb{P}(D_{I_*}\cdot \hat{v}=w)\leq\max_{w'}\mathbb{P}((D_{I_*})_{S^c\times S}\cdot \hat{v}_{S}=w')\leq \rho(\hat{v}_S)^{|I_*|-|S|}.
$$ Suppose that $|T|+K\leq|\hat{v}|$, then we have $\rho(\hat{v}_S)\leq c\rho_\mu(\hat{v}_S)\leq c\rho_\mu(u)$. Then applying the above two estimates to each $s=(T,u)$ we have
$$\begin{aligned}&\sum_{\hat{v}\in f^{-1}(s),|\hat{v}|\geq |T|+K}\max_{w\in\mathbb{F}_p^n}\mathbb{P}(D_{I_*}\cdot\hat{v}=w)\leq |f^{-1}(s)|\rho(u)^{|I_*|-|S|}\\&
\leq (\frac{2}{\rho_\mu(u)})^{|I_*|-|T|}(c\rho_\mu(u))^{|I_*|-|T|-K}\leq 2^{-4|I_*|}, 
\end{aligned}$$where we assume that $c<\frac{1}{32}$ and for the final inequality we use $\rho_\mu(u)\geq 1/p$ and thus $\rho_\mu(u)^K\geq 2^{-|I_*|}$. That is, we have 
\begin{equation}\label{caseB2}
\sum_{s\in\mathcal{S}}\sum_{\hat{v}\in f^{-1}(s):|\hat{v}|\geq |T|+K}\max_w\mathbb{P}(D_{I_*}\cdot \hat{v}=w)\leq 2^{-|I_*|}.
\end{equation}
On the other hand, when $|T|=t$ and $|\hat{v}|\leq t+K$, we use $\rho(\hat{v}_S)\leq \rho_\mu(\hat{v}_S)\leq\rho_\mu(u)\leq (1-\mu/2)^t$, where the last inequality follows from Lemma \ref{lemma3.97}. Then we have at most $p^{t+K}\binom{|I_*|}{\leq t+K}$ choices for $\hat{v}$. Then we have 
\begin{equation}\label{caseB3}\begin{aligned}&
\sum_{s\in\mathcal{S}}\sum_{\hat{v}\in f^{-1}(s),|\hat{v}|\leq |T|+K}\max_w\mathbb{P}(D_{I_*}\cdot\hat{v}=w)\\&\leq \sum_{t=1}^{\frac{2}{\nu}\log p} p^t\binom{n}{\leq t+K}(1-\mu/2)^{t(|I_*|-d-K)}\leq(1-\mu)^{|I_*/8|}.\end{aligned}
\end{equation}
\textbf{Step 3: Combining the estimates.} We have enumerated all candidates for the vector $v_{I_*}$, and we now find the cardinality for the possible vectors $v_{I_1},v_{I_2},\cdots,v_{I_{n_*-1}}$.

Recall that we have assumed $\rho_\mu(v_{I_k})\in[2^{-t_k},2^{-t_k+1}]$ for each $k$. Set $D_k=\frac{2}{\mu}(t_k+1)$ for each $k$, applying Theorem \ref{inversetheorem}, we get as in \eqref{balloonsinthesky} that the total cardinality for these vectors $v_{I_k}$ is at most: 
\begin{equation}\label{B4eq}\#\{\text{Candidates for } v_{I_k}\}\leq 
2^{|I_k|+9}p^{2D_k}(2^{\mu D_k/2})^{|I_k|-D_k}\leq 2^{|I_k|+9}\cdot2^{(t_k+1)|I_k|}\cdot e^{\frac{8}{\mu}(\log p)^2}.
\end{equation}

 Combining \eqref{caseB1}, \eqref{caseB2}, \eqref{caseB3} and the cardinality bound \eqref{B4eq},
we deduce that
\begin{equation}\label{BCthesecondline}\begin{aligned}&
\mathbb{P}(\exists v\in \ker S_n:\rho_\mu(v_{I_k})\in[2^{-t_k},2^{-t_k+1}]\forall k,|v_{I_k}|\geq K\forall k)\\&\leq c^{n-|I_*|}2^{20n}e^{\frac{2}{\mu}(\log p)^2(\frac{n}{e_n}-1)}\cdot (2^{-|I_*|}+(1-\mu)^{|I_*|/8}).\end{aligned}
\end{equation} By our assumption, $(\log p)^2\frac{n}{e_n}\leq \rho^2n^\alpha\frac{n}{n^\alpha/s}\leq\rho^2 sn$ for $s$ fixed (where $e_n=\lfloor d_n/s\rfloor$), so that by setting $c>0$ sufficiently small (which fixes a value of $\mu=\mu(c)$) and then setting $\rho>0$ sufficiently small relative to $\mu=\mu(c)$[more formally, we can find a function $\rho(\mu)$ such that wheenver $\rho\leq\rho(\mu)$ the following holds], we can guarantee that the second line of \eqref{BCthesecondline} is at most $\exp(-\Omega(n))$. Since the number of dyadic decompositions is at most $2^{o(n)}$ by \eqref{intervalsdyadicdecom}, we conclude that 
\begin{equation}\label{B613}
    \mathbb{P}(\exists v\in \ker S_n\text{ of type II}:|v_{I_k}|\geq K\forall k)\leq\exp(-\Omega(n)).
\end{equation}

 We now treat the leftover case where for some $k$, we have $v_{I_k}=0$ or $0<|v_{I_k}|<K$ for the constant $K>0$ in Lemma \ref{lemmafourier}. Denote by $\mathcal{R}$ the subset of indices $k$ such that $|v_{I_k}|\leq K$. Then 
 whenever $0<|v_{I_k}|\leq K$ we use the trivial fact that $\rho(v_{I_k})\leq\frac{1}{2}$ combined with the following upper bound on the number of such $v_{I_k}$: the cardinality is bounded by $$\#\{\text{Candidates for } v_{I_k}:|v_{I_k}|\leq K\}\leq(e_n)^Kp^K=e^{O_K(n^{\alpha/2}\log n)}.$$ If $v_{I_k}=0$, then we proceed trivially. The similarly as in estimate \eqref{B613}, we deduce that for any possible choice of the set $\mathcal{R}$, we have
 \begin{equation}\label{B61345}
    \mathbb{P}(\exists v\in \ker S_n\text{ of type II}:|v_{I_k}|\geq K\text{ if and only if } k\notin\mathcal{R})\leq\exp(-\Omega(n)).
\end{equation}
Enumerating all the possible subsets $\mathcal{R}\subset[\lfloor n/e_n\rfloor+1]$, we have considered all vectors of type II, and thus we complete the proof of the lemma. \end{proof}

\subsubsection{Eliminating kernel vectors of structure type III}

We need the following generalization of vectors of type II, which we define as type III below:

\begin{notation}\label{notationcases3}
     We say that a vector $v\in\mathbb{F}_p^n\setminus\{0\}$ is of type III$(\mu)$ if (i)
for at least a $\frac{1}{4}$ fraction of the indices $k$, $1\leq k\leq \lfloor n/e_n\rfloor+1$, we have $\rho_\mu(v_{I_k})\geq\frac{e_n}{p}$, and (ii), there exists two such consecutive indices $k$. That is, 
$$\{v\text{ is of type III$(\mu)$}\}\Leftrightarrow\begin{cases}\#\{k:\rho_\mu(v_{I_k})\geq \frac{e_n}{p}\}\geq \frac{1}{4}\frac{n}{e_n}+1,\\\exists k:\rho_\mu(v_{I_k})\geq\frac{e_n}{p},\rho_\mu(v_{I_{k+1}})\geq\frac{e_n}{p},\\\end{cases}$$where we identify $v_{I_{\lfloor n/e_n\rfloor+2}}$ with $v_{I_1}$.
\end{notation}
Then we can similarly prove the following lemma: 

\begin{lemma}\label{newtype3}We can find a constant $c_0>0$ such that for any $c\in(0,c_0)$, we can find $\mu=\mu(c)\in(0,\frac{1}{4}]$ such that, whenever $\rho\leq\rho(\mu)$ [where $\rho(\mu)$ is a fixed function of $\mu$ only] and $\tau=\rho/32$,
 then the following holds with probability $1-\exp(-\Omega(n))$: There exists no vector $v\in\mathbb{F}_p^n\setminus\{0\}$ of type III$(\mu)$, such that $S_n\cdot v=0$.     
\end{lemma}

\begin{proof}
Since the variance profile of $\mathcal{D}$ is translation invariant, we can assume without loss of generality that $\rho_\mu(v_{I_*})\geq\frac{e_n}{p}$, $\rho_\mu(v_{I_1})\geq\frac{e_n}{p}$ and  $\rho_\mu(v_{I_1})\geq \rho_\mu(v_{I_*})$, where $I_*=I_{\lfloor n/e_n\rfloor+1}$ is the last block.

Let $\Omega_1,\Omega_2,\Omega_3$ be a partition of $\lfloor n/e_n\rfloor$ such that $k\in \Omega_1$ if $\rho_\mu(v_{I_k})\geq\frac{e_n}{p}$, $k\in\Omega_2$ if $\rho_\mu(v_{I_k})<\frac{2}{p}$ and $k\in\Omega_3$ if $\frac{2}{p}\leq\rho_\mu(v_{I_k})<\frac{e_n}{p}$. Then $|\Omega_1|\geq\frac{1}{4}\frac{n}{e_n}$ since $v$ is of type III.

Then we can deduce as in \eqref{caseB1} that
\begin{equation}\label{section2caseB1}\begin{aligned}&
\mathbb{P}(S_n\cdot v=0)\leq \prod_{k=2}^{\lfloor n/e_n\rfloor+1}\rho(v_{I_k})^{|I_{k-1}|} \cdot \sup_w\mathbb{P}(D_{\lfloor n/e_n\rfloor+1}\cdot v_{I_{\lfloor n/e_n\rfloor+1}}=w)
 \\&\leq c^{\frac{1}{4}n}\cdot \prod_{k=2}^{\lfloor n/e_n\rfloor+1}\rho_\mu(v_{I_k})^{|I_{k-1}|} \cdot \sup_w\mathbb{P}(D_{\lfloor n/e_n\rfloor+1}\cdot v_{I_{\lfloor n/e_n\rfloor+1}}=w)\\&\leq c^{\frac{1}{4}n}\cdot \prod_{k=1}^{\lfloor n/e_n\rfloor}\rho_\mu(v_{I_k})^{|I_{k}|} \cdot \sup_w\mathbb{P}(D_{\lfloor n/e_n\rfloor+1}\cdot v_{I_{\lfloor n/e_n\rfloor+1}}=w),   \end{aligned} \end{equation}
where for any given $c>0$, if $k\in \Omega_1$ we apply Lemma \ref{lemmafourier} to find $\mu\in(0,\frac{1}{4}]$ such that $\rho(v_{I_k})\leq c\rho_\mu(v_{I_k})$, and if $k\in\Omega_2$ or $k\in\Omega_3$ we use the trivial bound $\rho(v_{I_k})\leq \rho_\mu(v_{I_k})$. The $c^{\frac{1}{4}n}$ term comes from $|\Omega_1||I_1|\geq\frac{1}{4}n$. In the last line we used $\rho_\mu(v_{I_1})\geq\rho_\mu(v_{I_*})$.

Meanwhile, as in \eqref{caseB2} and \eqref{caseB3}, we can deduce that
\begin{equation}
    \label{lines701}\sum_{\hat{v}\in\mathbb{F}_p^{|I_*|}:\rho_\mu(\hat{v})\geq \frac{e_n}{p}}\max_w\mathbb{P}(D_{I_*}\cdot \hat{v}=w)\leq 2^{-|I_*|}+(1-\mu)^{|I_*|/8},
\end{equation}which enumerates all the possible choices of $v_{I_*}$.

For each $k\in\Omega_1\cup\Omega_3$, we find an integer $t_k$ such that $\rho_\mu(v_{I_k})\in[2^{-t_k},2^{-t_k+1}]$ for each $k$. Set $D_k=\frac{2}{\mu}(t_k+1)$ for each $k$.
The cardinality of candidates of $v_{I_k}$ can be bounded as in \eqref{B4eq}. If $k\in\Omega_2$, then $v_{I_k}$ can range over all $p^{|I_k|}$ possible vectors. Then, for this choice of $\Omega_1,\Omega_2$ and $\{t_k\}_{k\in\Omega_1}$, we can combine this cardinality bound for $v_{I_k}$ with \eqref{section2caseB1} and \eqref{lines701} to deduce that, similar to \eqref{BCthesecondline},
\begin{equation}\label{newBCthesecondline}\begin{aligned}&
\mathbb{P}(\exists v\in \ker S_n:\rho_\mu(v_{I_k})\in[2^{-t_k},2^{-t_k+1}]\forall  k\in\Omega_1\cup\Omega_3,\rho_\mu(v_{I_k})\leq\frac{2}{p}\forall k\in\Omega_2)\\&\leq c^{\frac{1}{4}n}2^{20n}e^{\frac{8}{\mu}(\log p)^2(\frac{n}{e_n})}\cdot (2^{-|I_*|}+(1-\mu)^{|I_*|/8}).\end{aligned}
\end{equation}
 Now we set $c>0$ sufficiently small (which fixes a value of $\mu=\mu(c)$) and then set $\rho>0$ sufficiently small relative to $\mu=\mu(c)$[more formally, we can find a function $\rho(\mu)$ such that whenever $\rho\leq\rho(\mu)$ the following holds], we can guarantee that the second line of \eqref{newBCthesecondline} is at most $\exp(-\Omega(n))$. Since the number of possible choices of $t_i$ and the choice of $\Omega_1,\Omega_2,\Omega_3$ is at most $2^{o(n)}$ by \eqref{intervalsdyadicdecom}, we take a union bound over all $t_i,\Omega_1,\Omega_2,\Omega_3$ and conclude that 
\begin{equation}\label{newB613}
    \mathbb{P}(\exists v\in \ker S_n\text{ of type III$(\mu)$})\leq\exp(-\Omega(n)).
\end{equation}
 This completes the proof.
\end{proof}

\subsubsection{Eliminating kernel vectors of structure type IV}

The vectors $v$ in Theorem \ref{theoremsymmetriccase} that are not covered by type I, II and III can be summarized as of type IV, as defined below:

\begin{notation}
    We say a vector $v\in\mathbb{F}_p^n\setminus\{0\}$ is of type IV$(\mu)$ if 
    $$\{v\text{ is of type IV$(\mu)$}\}\Leftrightarrow\begin{cases}\exists k_1:\rho_\mu(v_{I_{k_1}})\leq\frac{e_n}{p},\rho_\mu(v_{I_{k_1+1}})\leq \frac{e_n}{p},\\\exists k_2:\rho_\mu(v_{I_{k_2}})\geq\exp(-\tau n^{\alpha/2}),\rho_\mu(v_{I_{k_2+1}})\geq\exp(-\tau n^{\alpha/2}).\\\end{cases}$$
\end{notation}

It is easy to check that Types I,II, III,IV defined above cover all the vectors that are eliminated in Theorem \ref{theoremsymmetriccase}:

\begin{fact}\label{fact3.12} Let $\mu\in(0,\frac{1}{4}]$ be any real number.
    For a vector $v\in\mathbb{F}_p^n\setminus\{0\}$ such that there exists some $k\in[\lfloor n/e_n\rfloor]$ with
$\min(\rho_\mu(v_{I_k}),\rho_\mu(v_{I_{k+1}}))\geq\exp(-\tau n^{\alpha/2}))$, then $v$ must be either of type I$(\mu)$ or of type II$(\mu)$ or of type III$(\mu)$ or of type IV$(\mu)$.
\end{fact}

\begin{proof}
    If $v$ is not of type III$(\mu)$ but satisfies the stated property, then $v$ satisfies the second condition in type III$(\mu)$ in Notation \ref{notationcases3}. Then $v$ violates the first condition in Notation \ref{notationcases3}, which means $$\#\{k:\rho_\mu(v_{I_k})\geq \frac{e_n}{p}\}< \frac{1}{4}\frac{n}{e_n}+1.$$Since there are $\lfloor \frac{n}{e_n}\rfloor +1$ total intervals, we must be able to find two consecutive intervals $v_{I_{k_1}}$, $v_{I_{k_1+1}}$ satisfying the first condition $\rho_\mu(v_{I_{k_1}})\leq\frac{e_n}{p},\rho_\mu(v_{I_{k_1+1}})\leq \frac{e_n}{p}$ of Type IV$(\mu)$. The second condition in Type IV$(\mu)$ is already assumed.
\end{proof}

\begin{lemma}\label{lemmainCase4}
We can find a constant $c_0>0$ such that for any $c\in(0,c_0)$, we can find $\mu=\mu(c)\in(0,\frac{1}{4}]$ such that, whenever $\rho\leq\rho(\mu)$ [where $\rho(\mu)$ is a fixed function of $\mu$ only] and $\tau=\rho/32$, the following holds with probability $1-\exp(-\Omega(n))$: There exists no vector $v\in\mathbb{F}_p^n\setminus\{0\}$ of type IV$(\mu)$, such that $S_n\cdot v=0$.

\end{lemma}

The proof of Theorem \ref{theoremsymmetriccase} is now immediate assuming this lemma:
\begin{proof}[\proofname\ of Theorem \ref{theoremsymmetriccase} assuming Lemma \ref{lemmainCase4}] We combine Lemma \ref{lemma3.9}, \ref{newtype2}, \ref{newtype3} and \ref{lemmainCase4} to deduce that, we can find a choice of $c_0\in(0,1)$ such that, for any $c\in(0,c_0)$ we can find $\mu=\mu(c_0)$ [This is the same function in all four lemmas] such that we can find a $\rho(\mu)>0$ [$\rho(\mu)$ is a function of $\mu$] satisfying that, for any $0<\rho<\rho(\mu)$ and $\tau=\rho/32$, the following holds with probability $1-\exp(-\Omega(n))$: there exists no vector $v\in\mathbb{F}_p^n\setminus\{0\}$ of either type I$(\mu)$ or II$(\mu)$ or III$(\mu)$ or IV$(\mu)$ such that $S_n\cdot v=0$.

    By Fact \ref{fact3.12}, these four types of vectors cover all vectors $v$ such that for some $k$, $\min(\rho_\mu(v_{I_k}),\rho_\mu(v_{I_{k+1}}))\geq\exp(-\tau n^{\alpha/2}))$. This is precisely the claim of Theorem \ref{theoremsymmetriccase}.
\end{proof}

Finally, we complete the proof of Lemma \ref{lemmainCase4}:

\begin{proof}[\proofname\ of Lemma \ref{lemmainCase4}]

Since the variance profile of $\mathcal{D}$ is translation invariant, we can assume without loss of generality that \begin{equation}\label{whatisrhomu1}\rho_\mu(v_1)\leq\frac{e_n}{p},\quad \rho_\mu(v_{I_*})\leq\frac{e_n}{p}\end{equation} (the first and last block) and that for some $\Delta\in[\lfloor n/e_n\rfloor]$ we have \begin{equation}\label{whatisrhomu2}\rho_\mu(v_{I_\Delta})\geq\exp(-\tau n^{\alpha/2}),\quad \rho_\mu(v_{I_{\Delta+1}})\geq\exp(-\tau n^{\alpha/2}).\end{equation}

For any two integers $x\leq y$, denote by $(D)_x^y:=\{D_x,D_{x+1},\cdots,D_y\}$ and $(U)_x^y:=\{U_x,U_{x+1},\cdots U_y\}$, then we apply an iterative conditioning argument as below to deduce that

Then we proceed with

 $$\begin{aligned}&
\mathbb{P}(S_n\cdot v=0)\leq\sup_w\mathbb{P}(\mathcal{D}\cdot v=w)
\\&\leq\prod_{k=\Delta+1}^{\lfloor n/e_n\rfloor+1}\sup_{w,(D)_1^{k-1},(U)_1^k}\mathbb{P}(U_{k}^Tv_{I_{k-1}}+D_{k}v_{I_k}+U_{{k+1}}v_{I_{k+1}}=w\mid (D)_1^{k-1},(U)_1^{k} )
\\&\cdot \sup_{w_1,\cdots,w_\Delta}\mathbb{P}(U_k^Tv_{I_{k-1}}+D_kv _{I_{k}}+U_{k+1}v_{I_{k+1}}=w_k\quad \forall 1\leq k\leq \Delta)
\\&\leq \left(\prod_{k=\Delta+2}^{\lfloor n/e_n\rfloor+1}\sup_w \mathbb{P}(U_k v_{I_k}=w)\right)\cdot\sup_w\mathbb{P}(D_{\lfloor n/e_n\rfloor+1}\cdot v_{I_{\lfloor n/e_n\rfloor+1}}=w)
\\&\cdot \prod_{k=1}^{\Delta}\sup_{w,(U)_{k+1}^{\Delta+1},(D)_{k+1}^\Delta}\mathbb{P}(U_k^Tv_{I_{k-1}}+D_kv_{I_k}+U_{k+1}v_{I_{k+1}}=w\mid (U)_{k+1}^{\Delta+1},(D)_{k+1}^{\Delta})\\&\leq \prod_{k=\Delta+2}^{\lfloor n/e_n\rfloor+1}\sup_w\mathbb{P}(U_kv_{I_k}=w)
\cdot\prod_{k=2}^\Delta\mathbb\sup_w\mathbb{P}(U_k^Tv_{I_{k-1}}=w)
\\&\cdot \sup_w \mathbb{P}(D_1v_{I_1}=w)\cdot \sup_w \mathbb{P}(D_{\lfloor n/e_n+1\rfloor}\cdot v_{I_{\lfloor n/e_n+1\rfloor}}=w).
    \end{aligned}, $$ 
Informally, after the above conditioning procedure, the vectors $v_{I_1},v_{I_{\lfloor n/e_n\rfloor+1}}$ appear twice but  $v_{I_\Delta}$ and $v_{I_{\Delta+1}}$ does not appear in the product: the former two vectors are well anti-concentrated but the latter two vectors are poorly anti-concentrated.

As proven in \eqref{howdowecanboundthis}, we have \begin{equation}\label{3.18}\sup_w\mathbb{P}(D_1v_{I_1}=w),
\sup_w\mathbb{P}(D_{\lfloor n/e_n\rfloor+1}\cdot v_{I_{\lfloor n/e_n\rfloor+1}}=w)\leq (\sqrt{\frac{2}{p}})^{|I_{\lfloor n/e_n\rfloor+1}|/2}.\end{equation} 
As proven in \eqref{cardinalityforv1}, we have 
\begin{equation}\label{cardinalityforv1v2v3}\#\{\text{Candidates for } v_{I_\Delta},v_{I_{\Delta+1}}\}\leq 2^{|I_1|}
p^{\frac{2\tau}{\mu} n^{\alpha/2}}\cdot p^{\frac{2\tau}{\mu} n^{\alpha/2}}\cdot (256\exp(\tau n^{\alpha/2}))^{|I_{1}|}.\end{equation}

We now enumerate the cardinality of the other components of $v_{I_k}$. Just as in the previous proofs, we need to find the cardinality of candidate vectors $v_{I_k}$ granted a knowledge of $\rho_\mu(v_{I_k})$. Let $\Omega_1,\Omega_2$ be a partition of $[\lfloor n/e_n\rfloor+1]$ such that $k\in\Omega_1$ if $\rho_\mu(v_{I_k})\geq\frac{2}{p}$ and $k\in\Omega_2$ if otherwise. 
For each $k\in\Omega_1$, we find an integer $t_k$ such that $\rho_\mu(v_{I_k})\in[2^{-t_k},2^{-t_k+1}]$ for each $k$. Set $D_k=\frac{2}{\mu}(t_k+1)$ for each $k$.
The cardinality of candidates of $v_{I_k}$ can be bounded as in \eqref{B4eq}. If $k\in\Omega_2$, then $v_{I_k}$ can range over all $p^{|I_k|}$ possible vectors. Then, for any choice of $\Omega_1,\Omega_2$ and $\{t_k\}_{k\in\Omega_1}$, we can combine this cardinality bound for $v_{I_k}$ with \eqref{3.18} and \eqref{cardinalityforv1v2v3} to deduce that, similar to \eqref{case1}, \begin{equation}\label{newcase1}\begin{aligned}&\mathbb{P}( \exists   v\in\ker S_n:\rho_\mu(v_{I_k})\in[2^{-t_k},2^{-t_k+1}]\forall k\in\Omega_1,\rho_\mu(v_{I_k})\leq\frac{e_n}{p}\forall k\in\Omega_2,\eqref{whatisrhomu1},\eqref{whatisrhomu2})\\&\leq 
2^n(\sqrt{\frac{2}{p}})^{2|I_1|/4} 2^{O_{\mu,\tau}(n^{\alpha/2}\log n)}\exp(2\tau n^{\alpha/2}|I_{1}
|).\end{aligned}\end{equation}
[Note: in this lemma we only use the basic inequality $\rho_\mu(v_{I_k})\leq \rho(v_{I_k})$ rather than the more refined version $\rho_\mu(v_{I_k})\leq c\rho(v_{I_k})$].
Here $a=O_{\mu,\tau}(b)$ means $a\leq Cb$ for $C$ a constant depending only on $\mu, \tau$ and $\alpha$. Recall that we take $p=\exp(\rho n^{\alpha/2})$ and $|I_{1}|\geq n^\alpha$ with $\alpha=\frac{2}{3}$. Then we see that as long as $\tau=\rho/32$, the above expression is super-exponentially small, which is actually at most $\exp(-\Omega(n^{1.5\alpha}))$. Finally, taking a union bound over all possible candidates of $t_i$ and $\Omega_1,\Omega_2$ completes the proof. 
\end{proof}

\subsection{Linear algebraic reductions and invertibility}\label{section3.223}
The proof of Theorem \ref{symmetricmaintheorem} relies on the following linear algebraic reductions:

\begin{fact}\label{fact3.2inclusions}(\cite{campos2021singularity}, Lemma A.6)
    If $\operatorname{rk}(S_n)=n-1$, then $\max_{i\in[n]}\operatorname{rk}(S_n^{(i)})\geq n-2$, where $S_n^{(i)}$ is the principal submatrix obtained from $S_n$ by removing its $i$-th row and column.
\end{fact}

\begin{fact}
    Since $S_n^{(i)}$ can be identified to have the same distribution as a version of $S_{n-1}$ (we still have $\{\pm 1\}$ random entries at given locations, while the value of the other fixed entries $c_{ij}$ may change) but with bandwidth $d_n-1$, and $d_n-1\geq \frac{1}{2}n^\alpha$ by definition, the structural result in Theorem \ref{theoremsymmetriccase} proven for $S_n$ can be immediately applied to $S_n^{(i)}$ with a minor change of constants, and with the following change of notations:

    We identify $[1,n]\setminus\{i\}$ with $[1,n-1]$ by replacing $x$ by $x-1$ if $x>i$,
    and then we partition $[1,n-1]=\cup_{k=1}^{\lfloor (n-1)/e_{n-1}\rfloor+1}I_k$, where the intervals $I_k$ are defined as in Theorem \ref{theoremsymmetriccase}. By definition of $e_n=\lfloor d_n/s\rfloor,s\geq 3$, we have $|I_k|\leq \frac{1}{2}d_n$ for each $k$, where $d_n$ is the bandwidth defined in Theorem \ref{symmetricmaintheorem}.

    We adopt these conventions implicitly throughout this section.  

\end{fact}

Then we decompose the event $\det S_n=0$ into the following cases:

\begin{lemma}\label{lemma3.16}
    For any symmetric random matrix $S_n$,
    $$\begin{aligned}&
\mathbb{P}(\det S_n=0)\leq \mathbb{P}(\operatorname{rk}(S_n)\leq n-2)\\&+\sum_{i=1}^n\mathbb{P}(\operatorname{rk}(S_n)=n-1,\operatorname{rk}(S_{n}^{(i)})=n-2)+\sum_{i=1}^n\mathbb{P}(\operatorname{rk}(S_n)=\operatorname{rk}(S_{n}^{(i)})=n-1).
\end{aligned}    $$
\end{lemma}
\begin{proof}
    This follows from the inclusion stated in Fact \ref{fact3.2inclusions}.\end{proof}

\begin{lemma}\label{lemma3.17}
    In the setting of Theorem \ref{symmetricmaintheorem}, we have 
    $$
\mathbb{P}(\operatorname{rk}(S_n)\leq n-2)\leq\exp(-\Omega(n^{\alpha/2})).
    $$
\end{lemma}

\begin{proof}
we reduce to the finite field $\mathbb{F}_p$.  If $\operatorname{rk}(S_n)\leq n-2$, then $\operatorname{rk}(S_n^{(i)})\leq n-2$ for each $i\in[n]$ and $S_n^{(i)}$ is not invertible. 
 Let $v\in\mathbb{F}_p^n\setminus\{0\}$ be such that $S_n\cdot v=0$, then $v$ must have a coordinate $i\in[n]$ such that $v_i\neq 0$. Without loss of generality we can assume $i=1$ and write $v=(v_1,w)$. Let $X$ be the first column of $S_n$ with the first entry removed, then we have $Xv_1+S_{n}^{(1)}\cdot w=0$. Let $u\in\mathbb{F}_p^n\setminus\{0\}$ be such that $S_n^{(1)}\cdot u=0$, then $u^T\cdot Xv_1+u^{T}S_n^{(1)}\cdot w=0$, so that $u^T\cdot X=0$ since $v_1\neq 0$. Then we can write 
 $$
\mathbb{P}(\operatorname{rk}(S_n)\leq n-2)\leq \sum_{i=1}^n\mathbb{P}(u^T\cdot X^i=0\text{ for any  $u\in\mathbb{F}_p^{n-1}\setminus\{0\}$ } :S_n^{(i)}\cdot u=0),
 $$where $X^i$ is the $i$-th row/ column of $S_n$ with the $i$-th entry removed. Consider the event \begin{equation}\label{satisfiesthe}\begin{aligned}\Omega_\tau:=&\{\forall 1\leq i\leq n,\forall u\in \ker S_n^{(i)}\setminus\{0\}:\\&\min(\rho_\mu(u_{I_k}),\rho_\mu(u_{I_{k+1}}))\leq\exp(-\tau n^{\alpha/2})\forall 1\leq k\leq \lfloor (n-1)/e_{n-1}\rfloor+1\}.\end{aligned}\end{equation}
    Then for a suitable choice of $\tau$, we have $\mathbb{P}(\Omega_\tau^c)\leq\exp(-\Omega(n))$ by Theorem \ref{theoremsymmetriccase}. On the event $\Omega_\tau$, we have 
    $$ u\text{ satisfies the second line of \eqref{satisfiesthe}
}\Rightarrow
\mathbb{P}_{X^i}(u^T\cdot X^i=0)\leq \exp(-\tau n^{\alpha/2}),
    $$ by definition of the random vector $X^i$. That is, 
    $$
\mathbb{P}(\operatorname{rk}(S_n)\leq n-2)\leq n\exp(-\tau n^{\alpha/2})+\exp(-\Omega(n))=\exp(-\Omega(n^{\alpha/2})).
    $$
\end{proof}

Then we upper bound the term $\mathbb{P}(\operatorname{rk}(S_n)=n-1,\operatorname{rk}(S_n^{(i)})=n-2)$. The following fact will be useful:

\begin{fact}\label{fact3,18}(\cite{campos2021singularity}, Lemma A.7) If $\operatorname{rk}(S_n^{(i)})=n-2$, then we can find a non-trivial column $a\in\mathbb{F}_p^{n-1}$ of $\operatorname{adj}(S_n^{(i)})$ so that $S_n^{(i)}\cdot a=0$ and if $\det(S_n)=0$, then $\sum_{j\in[n]\setminus\{i\}}a_jx_j=0$. Here $x=(x_1,\cdots,x_n)$ is the $i$-th row of $S_n$ and $a=(a_1,\cdots,a_n)$.
\end{fact}

\begin{lemma}\label{lemma3.19}
 In the setting of Theorem \ref{symmetricmaintheorem}, we have 
        $$
\sum_{i=1}^n\mathbb{P}(\operatorname{rk}(S_n)=n-1,\operatorname{rk}(S_n^{(i)})=n-2)\leq\exp(-\Omega(n^{\alpha/2})).
        $$
\end{lemma}

\begin{proof}
    By Fact \ref{fact3,18}, to bound the probability that $\operatorname{rk}(S_n)=n-1$ and $\operatorname{rk}(S_n^{(i)})=n-2$, we only need to bound the probability that there is a vector $a\in\mathbb{F}_p^{n-1}\setminus\{0\}$ such that $S_{n}^{(i)}\cdot a=0$ and $a\cdot X=0$, where $X$ is the $i$-th row of $S_n$ with the $i$-th entry removed and is thus independent of $S_n^{(i)}$.

    Recall the event $\Omega_\tau$ defined in \eqref{satisfiesthe}. Then on $\Omega_\tau$, for this fixed vector $a$ we have $\mathbb{P}_X(a\cdot X=0)\leq \exp(-\tau n^{\alpha/2})$. The vector $a$ is unique up to scalars since $\operatorname{rk}(S_n^{(i)})=n-2$, so summing the above probability with $\mathbb{P}(\Omega_\tau^c)\leq \exp(-\tau n^{\alpha/2})$ completes the proof. \end{proof}

Finally, we estimate the probability $\mathbb{P}(\operatorname{rk}(S_n)=\operatorname{rk}(S_n^{(i)})=n-1)$. We first need the following decoupling procedure:
 \begin{lemma}\label{lemma3.200}
For a vector $v\in\mathbb{F}_p^m$ and $J\subset[m]$, we denote by $v_J\in\mathbb{F}_p^{|J|}$ the restriction of $v$ onto the coordinates of $J$, and denote by $v_J^*$ the vector in $\mathbb{F}_p^m$ whose $i$-th coordinate is $v_i\cdot\mathbf{1}(i\in J)$. Consider a random vector $w\in\{-2,0,2\}^{n-1}$ with independent coordinates $(w_i)_{1\leq i\leq n-1}$ satisfying 
\begin{equation}
\label{definitionofw}\begin{cases}   
\mathbb{P}(w_i=2)=\mathbb{P}(w_i=-2)=\frac{1}{4},\quad \mathbb{P}(w_i=0)=\frac{1}{2},\quad \forall i\leq d_n\text{ or } i\geq n-d_n,\\
w_i=0\quad\forall d_n+1\leq i\leq n-d_n-1.
\end{cases}\end{equation} 
     Then, for any non-trivial partition $I,J$ with $I\cup J=[n-1]$, we have
     $$
\mathbb{P}(\operatorname{rk}(S_n)=\operatorname{rk}(S_n^{(i)})=n-1)\leq 2\cdot\mathbb{E}[\max_{a\in\mathbb{F}_p}\mathbb{P}(z_I\cdot w_I=a\mid S_n^{(i)})^{1/4}
\mathbf{1}[\operatorname{rk}(S_n^{(i)})=n-1]
],$$
where the expectation is over $S_n^{(i)}$ and $z:=(S_n^{(i)})^{-1}\cdot w_J^*$. 
 \end{lemma}   

    The proof of this lemma essentially follows from \cite{campos2021singularity}, Lemma A.9. But as the band matrix structure here is somewhat different, we write down a complete proof here. First we need the following decoupling lemma:

    \begin{lemma}\label{lemma3.21111}
        (Lemma 4.7 of \cite{costello2006random}) Let $X$ and $Y$ be independent random vectors, and $\mathcal{E}(X,Y)$ be an event which depends on $X$ and $Y$. Then let $X',Y'$ be independent copies of $X$ and $Y$, we have 
        \begin{equation}
\mathbb{P}(\mathcal{E}(X,Y)\leq \mathbb{P}(\mathcal{E}(X,Y)\cap\mathcal{E}(X,Y')\cap\mathcal{E}(X',Y)\cap\mathcal{E}(X',Y'))^\frac{1}{4}.            
        \end{equation}
    \end{lemma}

\begin{proof}[\proofname\ of Lemma \ref{lemma3.200}]
Without loss of generality, we will take $i=1$ when writing out the proof.
Recall that in the symmetric matrix $S_n$, we assumed $s_{ij}=u_{ij}+c_{ij}$ when $|i-j|\leq d_n$ or $|n-i+j|\leq d_n$ ,where $u_{ij}$ is uniform on $\{-1,1\}$ and $c_{ij}$ is independent of $u_{ij}$. Throughout the proof we first condition on the values of $c_{ij}$ and the values of $(s_{ij})_{d_n<|i-j|<n-d_n}$ and thus assume they are fixed throughout. That is, the $(1,1)$-th entry of $S_n$ is $c_{11}+u_{11}$ where now $c_{11}$ is a fixed constant. The entries of first row of $S_n$ are all conditioned to be fixed except the random variables $(u_{1j})_{1\leq j\leq d_n+1\text{ or } j\geq n-d_n+1}$ which are still random and uniform on $\{-1,1\}$.

Take $X=(q_i)_{i\in I}$ and $Y=(q_i)_{i\in J}$ be two families of random variables so that $q\in\mathbb{F}_p^{n-1}$ is determined by $X$ and $Y$. Then for each choice of $S_{n-1}^{(1)}$, we define the following event
$$
\mathcal{E}(X,Y):=\{\exists v\in\mathbb{F}_p^{n-1}:S_n^{(i)}\cdot v=q\text{ and }q\cdot v\in\{-1+c_{11},1+c_{11}\}\},
$$ which is determined by $X$ and $Y$.
We claim that if $\operatorname{rk}(S_n^{(1)})=n-1$ and the first row of $S_n$ is $(x_1+c_{11},q_1\cdots,q_{n-1})$ for some $x_1\in\{-1,1\}$ and some fixed $c_1$, then 
$$
\{\operatorname{rk}(S_n)=n-1\}\Rightarrow \{u\in\mathcal{E}(X,Y)\}.
$$
Actually, since $\det(S_n)=0\neq\det(S_n^{(1)})$, we can find a vector $v\in\mathbb{F}_p^n$ with $S_n\cdot v=0$ and $v_1=-1$. Denote by $v'=(v_2,\cdots,v_n)$, then $S_n^{(1)}\cdot v'=q$ and $q\cdot v'\in\{-1+c_{11},1+c_{11}\}$.

Now let $X',Y'$ be independent copies of $X$ and $Y$, we define the event
$$
\mathcal{E}_1:=\mathcal{E}(X,Y)\cap\mathcal{E}(X',Y)\cap\mathcal{E}(X,Y')\cap\mathcal{E}(X',Y').
$$ Applying Lemma \ref{lemma3.21111}, we deduce that 
$$
\mathbb{P}(\mathcal{E}(X,Y)\mid S_{n}^{(1)})\leq\mathbb{P}(\mathcal{E}_1\mid S_{n}^{(1)})^\frac{1}{4},
$$and thus 
$$
\mathbb{P}(\operatorname{rk}(S_n)=\operatorname{rk}(S_n^{(1)})=n-1)\leq \mathbb{E}[\mathbb{P}(\mathcal{E}_1\mid S_n^{(1)})^\frac{1}{4}\mathbf{1}[\operatorname{rk}(S_n^{(1)})=n-1]].
$$
Therefore, it suffices to prove that 
\begin{equation}
\label{equationin864}\mathbb{P}(\mathcal{E}_1\mid S_n^{(1)})\leq 16\cdot\max_{a\in\mathbb{F}_p}\mathbb{P}(z_I\cdot w_I=a\mid S_n^{(1)}),
\end{equation}
for any $S_n^{(1)}$ with $\operatorname{rk}(S_n^{(1)})=n-1$. To prove this, we take any such $S_n^{(1)}$ and denote by $Q:=(S_n^{(1)})^{-1}$ and $D:=\{-1+c_{11},1+c_{11}\}$. We show that if $q\in\mathcal{E}(X,Y)$ then $q^TQq\in D$. Indeed, we check that 
$$
q^{T}Qq=q^TQ\cdot S_n^{(1)}v'=q^Tv'\in\{-1+c_{11},1+c_{11}\}=D.
$$
As we used the notation $q=q(X,Y)$, we define $f(X,Y):=q^TQq$, then when the event $\mathcal{E}_1$ holds, we have 
$$
f(X,Y)-f(X',Y)-f(X,Y')+f(X',Y')\in 2D-2D.
$$
Let $\hat{w}:=q-q'$ where $q'$ is an independent copy of $q$, by a standard computation we verify that 
$$\begin{aligned}&
f(X,Y)-f(X',Y)-f(X,Y')+f(X',Y')\\&=2\sum_{i\in I}\sum_{j\in J}Q_{ij}(q_i-q_i')(q_j-q_j')=2z_I\cdot \hat{w}_I,\end{aligned}
$$where $z:=(S_n^{(i)})^{-1}\cdot \hat{w}_J^*$ so that $z_i=\sum_{j\in J}Q_{ij}\hat{w}_j=\sum_{j\in J} Q_{ij}(q_j-q_j')$.

Since $|D|=2$, this verifies \eqref{equationin864}. Recall that we assume the first row of $S_n$ is $(x_1,q_1,\cdots,q_{n-1})$ and we have conditioned on the entries of $S_n$ except the $u_{ij}$ variables such that $|i-j|\leq d_n$ or $|n-i+j|\leq d_n$, we deduce that $\hat{w}:=q-q'$ has the same law as the random vector $w$ in \eqref{definitionofw}. This completes the proof.
\end{proof}

We also use the following fact:
\begin{fact}\label{greenlandfact}[\cite{campos2021singularity}, Lemma A.10]
    For any $v\in\mathbb{F}_p^n$ and a partition $I\cup J=[n]$, we have
    $$
\rho_{1/2}(v)\leq\rho(v),\quad\text{and}\quad \rho(v)\leq \rho(v_I)\leq 2^{|J|}\rho(v).
    $$
\end{fact}

Finally, we are able to settle the last remaining case for $\det S_n=0$:
\begin{lemma}
In the setting of Theorem \ref{symmetricmaintheorem}, with $z_I,w_I$ defined as in Lemma \ref{lemma3.200} and intervals $I_1,I_2$ defined as in \eqref{911i1i2}, we have that 
\begin{equation}
 \label{kemma3.23greatworlds}\begin{aligned}   \mathbb{E}&[\max_{a\in\mathbb{F}_p}\mathbb{P}(z_I\cdot w_I=a\mid S_n^{(i)})^\frac{1}{4}\mathbf{1}[\operatorname{rk}(S_n^{(i)})=n-1]]\\&\leq (2^{|I_*|}\exp(-\tau n^{\alpha/2})+2^{-|J_*|})^\frac{1}{4}+3^{|J_*|}\exp(-\Omega(n^{\alpha/2})),
\end{aligned}\end{equation}where we denote by $$ I_*=(I_1\cup I_2)\setminus I,\quad J_*=J\setminus[d_n+1,n-d_n-1].
$$
\end{lemma}

\begin{proof}
In the setting of Lemma \ref{lemma3.200}, recall that $1\leq k\leq n-2$ and $J\subset[n-1]$ with $|J|=k$.  Recall that $S_n^{(i)}\cdot z=w_J^*\in W(J)$, where 
$$
W(J):=\{v\in\{-2,0,2\}^{n-1}:v_j=0\forall j\notin J,v_j=0\forall d_n+1\leq j\leq n-d_n-1 \}.
$$
We consider the following event, with $\beta=\exp(-\tau n^{\alpha/2})$,
$$\begin{aligned}
\mathcal{U}^{(J)}:=&\{\text{for every }v\in\mathbb{F}_p^{n-1}\setminus\{0\},S_n^{(i)}\cdot v\in W(J):\\&\min(\rho(v_{I_k}),\rho(v_{I_{k+1}}))\leq\beta\quad \forall 1\leq k\leq \lfloor (n-1)/e_{n-1}\rfloor+1\}.
\end{aligned}$$ Since $W(J)\leq 3^{|J\setminus[d_n+1,n-d_n-1]|}$, we use Theorem \ref{theoremsymmetriccase} to deduce that 
\begin{equation}\label{3casesA}\mathbb{P}(S_n^{(i)}\in \mathcal{U}^{(J)})\geq 1-3^{|J\setminus[d_n+1,n-d_n-1]|}\exp(-\Omega(n^{\alpha/2})).\end{equation}
We also check that
\begin{equation}
\label{3casesB}\operatorname{rk}(S_n^{(i)})=n-1\Rightarrow\mathbb{P}(z=0\mid S_n^{(i)})\leq 2^{-|J\setminus[d_n+1,n-d_n-1]|},
\end{equation}this is because if $z=0$ then $w_J^*=0$, which happens with probability $2^{-|J\setminus[d_n+1,n-d_n-1]|}$ because we have $w_i=0$ with probability $1/2$ for each $i\in J\setminus[d_n+1,n-d_n-1]$.

Recall that $w_J$ and $S_n^{(i)}$ together determine the vector $z$ and the entries of $w_I$ are independent of $w_J$. Recall that we consider the intervals $I_k$ in the decomposition $[n-1]=\cup_{k=1}^{\lfloor (n-1)/e_{n-1}\rfloor+1}I_k$, we have by definition that $\rho_{1/2}(z_{I\cap I_k})=\max_{a\in\mathbb{F}_p}\mathbb{P}(z_{I\cap I_k}\cdot w_{I\cap I_k}=a)$. Therefore, applying Fact \ref{greenlandfact}, we have for each $k$ that
\begin{equation}
\rho_{1/2}(z_{I\cap I_k})\leq \rho(z_{I\cap I_k})\leq 2^{|I_k\setminus I|}\rho(z_{I_k})
.\end{equation}
Since $S_n^{(i)}\in\mathcal{U}^{(J)}$ and $S_n^{(i)}\cdot z\in W(J)$, if $z\neq 0$ then we have $\min(\rho(z_{I_k}),\rho(z_{I_{k+1}}))\leq\beta$ for each $k$.
Then we can upper bound the quantity below, using randomness of $w_I$ at only $w_{I\cap I_1},w_{I\cap I_2}$ (since $|I_1|+|I_2|\leq d_n$): whenever $S_n^{(i)}\in\mathcal{U}^\beta$, the following estimate holds:
\begin{equation}
    \label{3caseC}\begin{aligned}
&\mathbb{P}(\{z_I\cdot w_I=a\}\cap\{z\neq 0\}\mid S_n^{(i)})\\&\leq \mathbb{E}[\min(\rho_{1/2}(z_{I\cap I_1}),\rho_{1/2}(z_{I\cap I_2}))\mathbf{1}(z\neq 0)\mid S_n^{(i)}]      
    \\&\end{aligned}\leq \max(2^{|I_1\setminus I|},2^{|I_2\setminus I|})\beta.
\end{equation}
Finally, the inequality \eqref{kemma3.23greatworlds} follows from combining \eqref{3casesA}, \eqref{3casesB} and \eqref{3caseC}.
\end{proof}

\begin{corollary}\label{corollary3.24}In the setting of Theorem \ref{symmetricmaintheorem}, we have
$$
\mathbb{P}(\operatorname{rk}(S_n)=\operatorname{rk}(S_n^{(i)})=n-1)=\exp(-\Omega(n^{\alpha/2}).
$$
\end{corollary}

\begin{proof}
    By Lemma \ref{lemma3.200}, we only need to show the right hand side of \eqref{kemma3.23greatworlds} is $\exp(-\Omega(n^{\alpha/2}))$ for some choice of partitions $I\cup J=[n-1]$. We can achieve this by taking $I\subset I_1\cup I_2$ (so that $I_*=0$) and take $|J_*| =\exp(-\tau'n^{\alpha/2})$ for some $\tau'>0$. We can verify that as long as $\tau'>0$ is sufficiently small, the right hand side of \eqref{kemma3.23greatworlds} is $\exp(-\Omega(n^{\alpha/2}))$. 
\end{proof}

Now we have completed the proof of Theorem \ref{symmetricmaintheorem}:

\begin{proof}
    [\proofname\ of Theorem \ref{symmetricmaintheorem}] This follows from combining Lemma \ref{lemma3.16}, \ref{lemma3.17}, \ref{lemma3.19} and Corollary \ref{corollary3.24}.
\end{proof}

\section*{Funding}
The author is supported by a Simons Foundation Grant (601948, DJ).

\printbibliography

\end{document}